\documentclass[a4paper,11pt]{amsart}
\usepackage{graphicx} % Required for inserting images
\usepackage{amsmath, amssymb, amsthm, float, amsfonts,url,faktor,calrsfs,mathtools} %making math pretty
\usepackage{enumerate}%better lists
\usepackage{extarrows,tikz-cd, xcolor,tikz} %diagrams
\usetikzlibrary{decorations.pathmorphing,graphs}
\usepackage{subfiles}%bookkeeping subfiles
\usepackage{caption}%allows captions to work better
\usepackage[labelformat=simple]{subcaption}
\usepackage{chngcntr} %allows xwithin commands to change counters on stuff %better formatting on a4
\usepackage{hyperref} %clickable everything
\DeclareMathAlphabet{\oldcal}{OMS}{zplm}{m}{n} %oldmathcal
\numberwithin{equation}{section} %counts eqns and figures within section.
\counterwithout{figure}{section}

\newcommand{\Z}{\mathbb{Z}}

\newcommand{\C}{\mathcal{C}}
\newcommand{\D}{\mathcal{D}}

\newcommand{\Mm}{\mathcal{M}}
\newcommand{\M}{\oldcal{M}}
\newcommand{\Sp}{\mathcal{S}}

\newcommand{\morC}[3][]{\oldcal{M}_{#1}(#2,#3)}
\newcommand{\norm}[1]{\vert #1 \vert}
\newcommand{\Cprod}{\mathcal{C}\otimes \mathcal{D}}
\newcommand{\Ztwo}{\Z/2\Z}

\newcommand{\noop}[1]{} 
\DeclareMathOperator{\Ob}{Ob}
\DeclareMathOperator{\Sq}{Sq}
\DeclareMathOperator{\sq}{sq}

\DeclareMathOperator{\Hom}{Hom}
\DeclareMathOperator{\gr}{gr}
\DeclareMathOperator{\source}{source}
\DeclareMathOperator{\target}{target}

\theoremstyle{plain}
\newtheorem{theorem}{Theorem}[section]
\newtheorem*{theorem*}{Theorem}

\newtheorem{lemma}[theorem]{Lemma}
\newtheorem{proposition}[subsubsection]{Proposition}

\theoremstyle{definition}
\newtheorem{definition}[theorem]{Definition}
\newtheorem{example}[theorem]{Example}

\theoremstyle{remark}
\newtheorem*{remark}{Remark}

\author[]{Chris Johnson}
\title{Products and Duality for Framed 1-flow Categories}
\date{\today}

\begin{document}
\begin{abstract}
   We describe a product between framed 1-flow categories and show that it satisfies the expected Cartan formula on its second Steenrod squares. We also describe a dual framed 1-flow category and show that, like the Spanier-Whitehead dual of spectra, the homology groups of a framed 1-flow category can be shown to be isomorphic to the cohomology groups of its dual and that this isomorphism can be shown to cooperate with the framed 1-flow category structure in a suitable sense.
\end{abstract}
\maketitle
\section{Introduction}

Framed 1-flow categories are a combinatorial simplification of the much more general, and earlier introduced \textit{framed flow categories} of Cohen, Jones and Segal (CJS) \cite{cjs}. These are graded categories in which morphisms between objects with gradings $n>m$ are described by a framed $(|n|-|m|-1)$-dimensional manifold with corners. Cohen, Jones and Segal show that one can use this data to produce a CW-complex whose homotopy type is unique up to stabilisation. This construction was notably used in knot theory to upgrade the Khovanov chain complex to a stable homotopy type for Khovanov homology \cite{Lipshitz_2014}, providing for each quantum grading a topological space whose cohomology is isomorphic to $Kh^{*,q}(L)$, and hence further refining the Jones polynomial to the Khovanov stable homotopy type. The notion of a 1-flow category was later introduced in the work of Lobb, Orson and Sch\"utz \cite{d33} where a combinatorial second Steenrod square was described and shown to be equivalent to the second Steenrod square arising from the CJS construction. This method of taking a chain complex, building a suitable framed flow category, and then dealing with the resulting suspension spectra from the CJS construction has found use across knot theory, in particular, see \cite{ bhattacharyya2025lipshitzsarkarstablehomotopytype, tao2026steenrodsquarelinkfloer, manolescu2021knotfloerstablehomotopy,rajapakse2026steenrodsquaresoddkhovanov}, and use of framed 1-flow categories was employed by Sch\"utz in \cite{schuetz2022secondsteenrodsquaresodd} to produce two valid cohomology operations on odd Khovanov homology, which more recently in \cite{rajapakse2026steenrodsquaresoddkhovanov} was shown to be equivalent to the second Steenrod squares on two distinct applications of the CJS construction to odd Khovanov homology. 

One thing missing from the picture is the study of framed flow categories themselves. In Khovanov homology, given a pair of links $L_1, L_2$, we have an isomorphism of the Khovanov chain complexes $CKh^{i,j}(L_1) \otimes CKh^{i,j}(L_2)$ with $CKh^{i,j}(L_1\sqcup L_2)$, and Lawson, Lipshitz and Sarkar showed in \cite{Lawson_2020} that this moreover corresponds to the smash product on the level of spectra, that is, $\oldcal{X}^q_{Kh}(L_1)\wedge\oldcal{X}^q_{Kh}(L_2) \cong \oldcal{X}^q_{Kh}(L_1\sqcup L_2)$. This was proven without the use of the framed flow category machinery, and so passing from the tensor product of the homology groups to the smash product of spectra did not necessitate the construction of an equivalent product of framed flow categories. Likewise, Khovanov homology and its stable homotopy type plays well with mirror knots too -- the cohomology theory flips in both degrees, and the stable homotopy invariant becomes its Spanier-Whitehead dual. 

Between the tensor product of chain complexes and smash product of spectra one might expect a similar product between framed flow categories. In Section \ref{sec:products} we define a product for the framed 1-flow categories and we justify its comparison with the smash product of spectra by showing that it behaves well with the combinatorial second Steenrod square, and we investigate its properties.

\begin{theorem*}
    Given framed 1-flow categories $\C$, $\D$, there exists a framed 1-flow category $\C\otimes \D$ for which the Cartan formula for the Steenrod operations is respected.
\end{theorem*}

One might also ask for a dual framed 1-flow category, in analogue with the Spanier Whitehead (SW) dual of spectra. For a spectrum $\oldcal{X}$, the SW dual $\oldcal{X}^\vee$ is characterised by the existence of a map $\mu: \oldcal{X}\wedge \oldcal{X^\vee} \rightarrow \oldcal{S}$, where $\oldcal{S}$ is the sphere spectrum, which induces an isomorphism between the homology of $\oldcal{X}$, and the cohomology of $\oldcal{X}^\vee$. 
\begin{theorem*}
    Let $\C$ be a framed 1-flow category. There exists a framed 1-flow category $\C^\vee$, and a framed mapping cone 1-flow category $\M_\mu$ which induces an isomorphism between $H_k(\C;\Z)$ and $H^{-k}(\C^\vee;\Z)$. Moreover, Steenrod squares are respected in a suitable sense.
\end{theorem*} 

\textbf{Acknowledgements:} The author would like to thank Dirk Sch\"utz for supervising this research and for providing many hours of support in the direction and proofs herein.

\textbf{AI declaration:} From proof-writing to proofreading, no AI was used in this paper at any stage of its construction. 

\section{Framed 1-flow Categories}
We recount in this section the definitions with mainly differences to the presentations in \cite{d33, Lipshitz_2014b}, and prove towards the end a lemma which we use in Section 3.

\begin{definition}\label{def:1flow}
A \textit{1-flow category} $\C$ consists of a finite set of objects $\Ob(\C)$, which come with a grading $|\cdot| : \Ob(\C) \to \Z$. The morphisms, which we call \textit{moduli spaces}, from any object $x$ to $y$ with $|x|\leq |y|+2$ are described by compact manifolds of dimension $|x|-|y|-1$. 
The 1-dimensional moduli spaces must satisfy a boundary condition:
If $a,c\in \Ob(\C)$ with $|a|-|c|=2$, then the boundary of the moduli space is given by
\[
 \partial \M(a,c)=\coprod_{b\in \Ob(\C),|b|=|a|-1} \M(b,c)\times \M(a,b)
\]
Moreover, while we do not have any 2-dimensional moduli spaces in a 1-flow category, the boundary condition inherited from the general flow categories on the 2-dimensional moduli spaces must still hold. That is, for $|a|-|d| =3$, we have: 
\begin{align*}
 \partial \M(a,d)= & \coprod_{b\in \Ob(\C),|b|=|a|-1} \M(b,d)\times \M(a,b) \\
 & \cup \coprod_{c\in \Ob(\C),|c|=|a|-2}\M(c,d)\times \M(a,c)
\end{align*}
where we note the pair of summed products share a common subset
\[
 \coprod_{\stackrel{(b,c)\in \Ob(\C)\times \Ob(\C)}{|b|=|c|+1=|d|+2}} \M(c,d)\times \M(b,c)\times \M(a,b),
\]
 which may be referred to as the \textit{vertices} of this boundary, alluding to the fact that these boundaries are all $2n$-gons (inc. 0-gons, or circles).
\end{definition}

Notably, 1-flow categories are not categories in the typical sense as they fail to respect composition between all objects. One could keep track of the 1-skeleton between all moduli spaces, but this does not have a clear benefit and so we stick with the definition proposed in \cite{d33}.

\begin{definition}
\label{def:signassignment}
 A \textit{sign assignment} on a 1-flow category is a function on the 0-dimensional moduli spaces which assigns to every $P$ a value $s(P)\in \mathbb{F}_2$ which must satisfy a compatibility condition on the boundaries of 1-dimensional moduli spaces, that is to say, if we have an interval component $I \subseteq \M(a,c)$, with boundary $\partial I = (P_1,Q_1) \sqcup (P_2,Q_2)$, then we require that 
 \[
 s(P_1) + s(P_2) + s(Q_1) + s(Q_2) = 1
 \]
\end{definition}

\begin{definition}
    For every 1-flow category $\C$ with a sign assignment, one can define an associated integral chain complex $C_*(\C;\Z)$. We have $C_k(\C;\Z) = \langle c\in\Ob(\C)\;|\;|c_i|=k\rangle$, with differential
    \[
    \partial(c) = \sum_{\substack{c'\in \Ob(\C)\\|c| = |c'|+1}} \sum_{P \in \M(c,c')} (-1)^{s(P)} c'.
    \]
    Likewise, we can define a cochain complex by defining generators for every $c\in \Ob(\C)$ $\phi_c:C_k(\C;\Z) \to \Z$ which are defined to evaluate as zero everywhere except at $\phi_c(c) = 1$ with the usual coboundary $\delta(\phi_c)[c'] = \phi_c(d(c'))$.
    
    We abuse notation, and will occasionally refer to objects of $\C$ as the generators of the cochain complex where it is clear to do so. 
\end{definition}

\begin{definition}
    For $P\in \M(x,y)$, we define $\source(P)= x$ and $\target(P) = y$, and $\gr_s(P)=|x|$ to be the grading of the `source' of $P$, and likewise $\gr_t(P) = |y|$ for the grading of the `target'. 
\end{definition}
\begin{definition}\label{framingdef}
    A \textit{pre-framing} $f$ on a 1-flow category is a function on the 1-dimensional moduli spaces which assigns to every component a value in $\mathbb{F}_2$.

   For $|a|=|d|+3$, let $f$ a pre-framing and $s$ a sign assignment for $\C$. Consider $\partial \M(a,d)$, a disjoint union of components $C_i$ which are either circles of the form $S \times \{R\}$ or $\{P\}\times T$, or $2n$-gons $K$. Note, every $2n$-gon is comprised of two `types' of intervals, those of the form $I\times \{R\}$, and $\{P\}\times J$. Suppose then that $f$ satisfies the following relation: 
    \[
    \sum_{C_i\in \partial\M(a,d)} 1+ \widetilde{f}(C_i) = 0,
    \]
    where $\widetilde{f}$ is defined to act on the connected components of $\partial\M(a,d)$ in the following fashion:
   $$\widetilde{f}(S\times \{Q\}) = f(S),$$
$$\widetilde{f}(\{P\}\times T) = f(T),$$
 $$ \widetilde{f}(K)= \sum_{I\times \{R\}\subset K}f(I) + \sum_{\{P\}\times J\subset K} 1 + s(P) + f(J),$$
    then we call $f$ a \textit{framing} for $\C$ and we call the triple $(\C,s,f)$ a \textit{framed 1-flow category}. 
\end{definition}
\begin{remark}
    The boundary condition effectively ensures that the boundary framed 1-manifolds are framed null cobordant. We follow the convention in place in \cite{schuetz2022secondsteenrodsquaresodd}, where 0-framed circles are counted as non-trivial- a slightly confusing quirk, owed to the fact that in some sense the tangent bundle is included in the `framing' of the moduli spaces.
\end{remark}
\begin{definition}
    A \textit{full subcategory} of a 1-flow category $(\Mm,s,f)$ is a 1-flow category $(\C,s_\C,f_\C)$ with $\Ob (\C)\subset \Ob (\Mm)$, and for any $x,y \in \Ob(\C)$, we have $\M_\Mm (x,y) = \M_{\C}(x,y)$. Moreover $\C$ must take the obvious restrictions of $s$ and $f$.

    If we have that for all $x \in \Ob(\C)$, any such $y$ with $\M_\Mm(x,y)\neq\emptyset$ is also in $\Ob(\C)$, we say $\C$ is a \textit{downwards-closed} subcategory of $\Mm$. Flipping the positions of $x$ and $y$, we call this an \textit{upwards-closed} subcategory of $\C$.   

     An \textit{inclusion} of 1-flow categories $i:\C\to \Mm$ is an identification of $\C$ with a downwards closed subcategory of $\Mm$, and likewise the \textit{projection} $:\Mm\to \D$ is induced by an the identification of the upwards-closed subcategory with $\D$ as a full subcategory of $\Mm$.

    For a triple $\C, \D, \Mm$ where $\C$ is a downwards closed subcategory of $\Mm$ and $\D$ is the \textit{complementary} upwards flow category, we write:
    \[
        0\longrightarrow \D \longrightarrow \Mm \longrightarrow \C \longrightarrow 0
    \]
    and call this a \textit{short exact sequence of 1-flow categories}. 
\end{definition} 
\begin{definition}
    Let $\varphi = \{x_i\}_{i=1}^n$ be a collection of objects of the same grading. Then for any $x\in \Ob(\C)$, $\M(x,\varphi)$ is defined as the disjoint union of (framed) moduli spaces $\M(x,x_i)$. Likewise for $\M(\varphi,y)$. 
\end{definition}
\subsection{Second Steenrod Squares}\label{subsec:secgraphstr}

We introduce the second Steenrod square mostly in line with Lobb, Orson, Sch\"utz \cite{d33} and Lipshitz, Sarkar \cite{Lipshitz_2014b}. 
The combinatorial second Steenrod square is defined using a gadget called a `special graph structure'. We will include their definition as in \cite{d33}.

\begin{definition}\label{def:graphstructure}
 A \em special graph structure \em $\Gamma=\Gamma(V,E,E',E'',L,s,f)$ consists of a set of vertices $V$, a function $s\colon V \to \{0,1\}$, a set of edges $E$, a subset $E'\subset E$ together with a function $f\colon E'\to \{0,1\}$, a subset $E'' \subset E-E'$ of directed edges, and a set of loops $L$. All sets are finite. It has to satisfy the following conditions:
\begin{enumerate}
 \item \label{item:gs1}Each vertex is contained in at least one edge, and at most in two edges. We denote the set of vertices contained in only one edge by $\partial V$ and call these the boundary points.
 \item \label{item:gs2}If $v\in \partial V$, then the unique edge $e(v)$ with $v\in e(v)$ satisfies $e(v)\in E'$.
 \item \label{item:gs3}If $v \in V-\partial V$, then there is an edge $e_1\in E'$ with $v\in e_1$, and an edge $e_2\in E-E'$ with $v\in e_2$.
 \item \label{item:gs4}If $e\in E'$ and $e=\{v_1,v_2\}$, then $s(v_1)\not=s(v_2)$.
 \item \label{item:gs5}If $e\in E-E'$ and $e=\{v_1,v_2\}$, then $s(v_1)=s(v_2)$ if and only if $e\in E''$.
\end{enumerate}
The edges $E$ determine an equivalence relation on the vertex set $V$ in the obvious way. The set of equivalence classes under this relation, taken together with the set $L$ of loops, forms the set of \em components \em of $\Gamma$. Notice that $L$ is a set independent from possible loops made from edges in $E$.
\end{definition}

\begin{definition}
    Let $\varphi$ be a cocycle of degree $k$ in a 1-flow category $\C$. Then there exists a pairing of the 0-dimensional moduli spaces $\M(x,\varphi)$, for each $x$ with $|x|=k+ 1$ (by definition). We call any such pairing an \textit{combinatorial matching} on the 0-dimensional moduli spaces $\M(x,\varphi)$, and for $P,Q$ in $\M(x,\varphi)$ matched we write $\overrightarrow{(P,Q)}$. We use $M_\varphi$ to refer to the set of all matched pairs associated to $\varphi$, and $M_\varphi^x$ to be the set of pairs $(P,Q)$ with $\source(P)=\source(Q)=x$. 
\end{definition}

\begin{example}\label{eg:first}
Let $\varphi$ be a cocycle of degree $k$ in $\C$ and choose a matching $M_\varphi$ on the 0-dimensional moduli spaces $\M(b,\varphi)$. Then, for $|a|=k+2$, we define a special graph structure $\Gamma (a,\varphi)$. 
\begin{enumerate}
    \item The vertices $V$ of $\Gamma$ are given by all pairs $\M(b,\varphi)\times\M(a,b)$, for all $b$ in degree $k+1$. The signs of these vertices is defined so $s((P,Q)) = s(P) + s(Q)$.
    \item The edges $E'$ are comprised of the intervals in $\M(a,\varphi)$. They are framed in $\Gamma$ by their framing in $\C$.
    \item The loop set $L$ is defined by the set of non-trivially framed circles in $\M(a,\varphi)$.
    \item The other undirected edges $E-E'$ are defined by an ordered matching of $\varphi$: if $\overrightarrow{(P,Q)}$ is one such matched pair in $\M(b,\varphi)$ and $s(P)\neq s(Q)$, then for any $R\subseteq \M(a,b)$ we have an interval $I\subset E'$ with boundary $(P,R)\sqcup (Q,R)$.
    \item  If $s(P)=s(Q)$, we instead add an interval as above to the set of directed edges $E''$ and we give the direction by the ordering of the matching $\overrightarrow{(P,Q)}$. 
\end{enumerate} 
\end{example}

\begin{definition}
Let $\Gamma=\Gamma(V,E,E',E'',L,s,f)$ be a special graph structure. For a component $C$ of $\Gamma$ which contains a vertex, let $F(C)$ be the sum of the framing values $f(e')$ where $e'\in E'$ is in $C$. Also, let $D(C)\in \Z/2$ be the number of oriented edges in $C$ which point in a chosen given direction.
\end{definition}

\begin{definition}\label{def:sqca}Let $(\C,s,f)$ be a framed $1$-flow category and let $\varphi$ be a cocycle in $C^k(\C;\Z/2\Z)$. Then define a cochain $\sq^\varphi\colon C_{k+2}(\C;\Z/2\Z)\to \Z/2\Z$ by 
\[
\sq^2_\varphi(a):=|L|+\sum_{C}\left(1+F(C)+D(C)\right)\in\Z/2,
\]
where $L$ is the loop set of $\Gamma_\mathcal{C}(a,\varphi)$ and the sum is taken over all components $C$ of~$\Gamma_\mathcal{C}(a,\varphi)$ containing a vertex.
\end{definition}

It was shown in \cite{d33} that this combinatorial second Steenrod square is equivalent to the second Steenrod square in \cite{Lipshitz_2014}.

Using the matching, we may also define a combinatorial $\sq^1$ operating on cochains in $C^k(\C;\Z)$. See, for instance \cite{Lipshitz_2014b} where it is also described.
\begin{definition}\label{def:sq1}
    Associated to the Bockstein homomorphism of the short exact sequence
    \[
    0\rightarrow \Z/2\Z\rightarrow \Z/4\Z\rightarrow \Z/2\Z\rightarrow0,
    \]
    we have, for $\varphi\in C^k(\C;\mathbb{F}_2)$, an operation $\sq^1_\varphi\in C^{k+1}(\C;\Z/2\Z)$ defined by its output on objects $x$ in degree $k+1$ by
      \[
        \sq^1_\varphi(x) = \sum_{\substack{(P,Q)\in M^x_\varphi\\ s(P)=s(Q)}} 1.
      \]

    One may also define a related operation $\widetilde{\sq}^1_\varphi \in \C^{k+1}(\C;\Z)$, which is associated to the exact sequence
    \[
    0\rightarrow \Z\rightarrow \Z\rightarrow \Z/2\Z\rightarrow0
    \]
    by
    \[
        \sq^1_\varphi(x) = \sum_{\substack{(P,Q)\in M\\ s(P)=s(Q)}} (-1)^{s(P)}.
    \]

    These two are related by change of coefficients (given by $p$ below), that is the following commutes for all $p$ and $\varphi$:
    \begin{center}
        \begin{tikzpicture}[font=\small]
            \node (a) at (0,0) {$H^k(\C;\Z/2\Z)$} ;
            \node (b) at (3,0) {$H^{k+1}(\C; \Z)$} edge [<-, "$\widetilde{\sq}^1_\varphi$"] (a);
            \node (c) at (3,-2) {$H^{k+1}(\C;\Z/2\Z)$} edge [<-, "$p$"] (b)
                                                        edge [<- , "$\sq^1_{\varphi}$"] (a);
        \end{tikzpicture}.
    \end{center}

\end{definition}

\begin{example} \label{ex:sq2cat}
    Let $\C$ be the 1-flow category with $\Ob(\C) = \{c_n, c_{n+1}^0,\linebreak c_{n+1}^1, c_{n+2}\}$, so that these objects are graded according to their subscripts. Define the 0-dimensional moduli spaces by the following diagram:
    
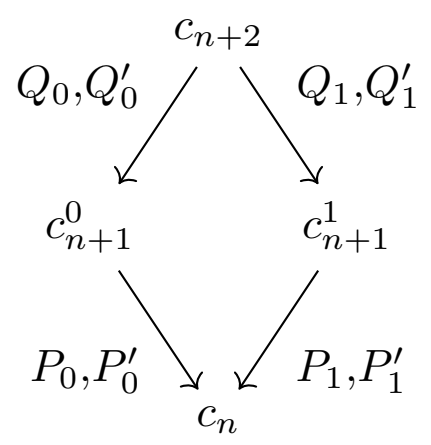
\begin{center}
            \begin{tikzpicture}[font =\small, transform shape,scale=0.85]
        \node (u) at (-3,1.5) {$c_n$};
                \node (v) at (-2,3) {$ c_{n+1}^1$} edge [->, "$P_1{,}P_1'$"] (u);
                \node (w) at (-4,3) {$ c_{n+1}^0$} edge [->, "$P_0{,}P_0'$", swap] (u);
                \node (x) at (-3, 4.5) {$c_{n+2}$} edge [-> , "$Q_1 {,} Q_1'$"] (v)
                                                         edge [->, "$ Q_0 {,} Q_0'$", swap] (w);
    \end{tikzpicture}
    \end{center}
    with all sign assignments equal to $0$ except $s(Q_0) = s(Q_0') =1$. We have then four intervals, as there are 8 boundary points in $\partial \M(c_{n+2},c_n)$. We define these intervals $I_1,I_2,I_3, I_4 \subseteq \M(c_{n+2},c_n)$ by their boundaries,
    \begin{align*}
        \partial I_1 &=  (P_0, Q_1) \sqcup (Q_0, P_1) 
        &\partial I_2 &=  (P_0, Q_1') \sqcup (Q_0', P_1)  \\
        \partial I_3 &=  (P_0', Q_1) \sqcup (Q_0, P_1') 
        &\partial I_4 &=   (P_0', Q_1') \sqcup (Q_0', P_1') \\
    \end{align*}%
    and set their framing identically equal to 0. So, we have
    
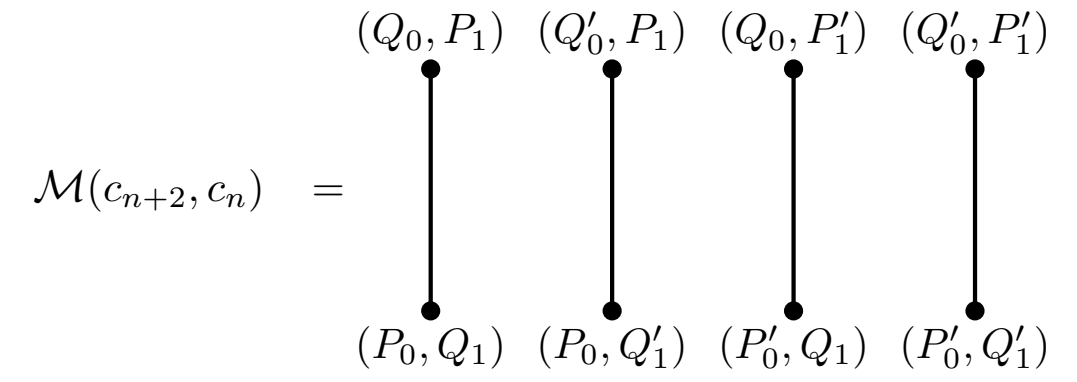
\begin{center}
        \begin{tikzpicture}[font=\small, scale =0.8, transform shape]
            \node (a) at (-3,0) {$\M(c_{n+2},c_n) \;\;\; =$};
            \filldraw[black] (-1,-1) circle (2pt);
            \node (b) at (-1,-1.3) {$(P_0,Q_1)$};
            \filldraw[black] (-1,1) circle (2pt);
            \node (c) at (-1,1.3) {$(Q_0, P_1)$} edge[-, thick] (b);
            \filldraw[black] (0.5,-1) circle (2pt);
            \node (b2) at (0.5,-1.3) {$(P_0,Q_1')$};
            \filldraw[black] (0.5,1) circle (2pt);
            \node (c2) at (0.5,1.3) {$(Q_0', P_1)$} edge[-, thick] (b2);
            \filldraw[black] (2,-1) circle (2pt);
            \node (b3) at (2,-1.3) {$(P_0',Q_1)$};
            \filldraw[black] (2,1) circle (2pt);
            \node (c3) at (2,1.3) {$(Q_0, P_1')$} edge[-, thick] (b3);
            \filldraw[black] (3.5,-1) circle (2pt);
            \node (b4) at (3.5,-1.3) {$(P_0',Q_1')$};
            \filldraw[black] (3.5,1) circle (2pt);
            \node (c4) at (3.5,1.3) {$(Q_0', P_1')$} edge[-, thick] (b4);
        \end{tikzpicture}
    \end{center}
        Now $(\C,s,f)$ satisfy the criteria to be a framed 1-flow category. We shall demonstrate how to compute the second Steenrod square in $\C$. There is only one graph structure to compute, as we have a unique object in degrees $n$, $n+2$. Firstly we need a combinatorial boundary matching, we set $M_{c_n} = \{(P_0,P_0'), (Q_0,Q_0')\}$.

        The graph $\Gamma_{c_n}(c_{n+2})$ is described by the following diagram, where dashed edges represent the edges in $(E-E') \cup E''$, however, since the signs are the same for the matched pairs, all dashed edges will in this case be in $E''$. %
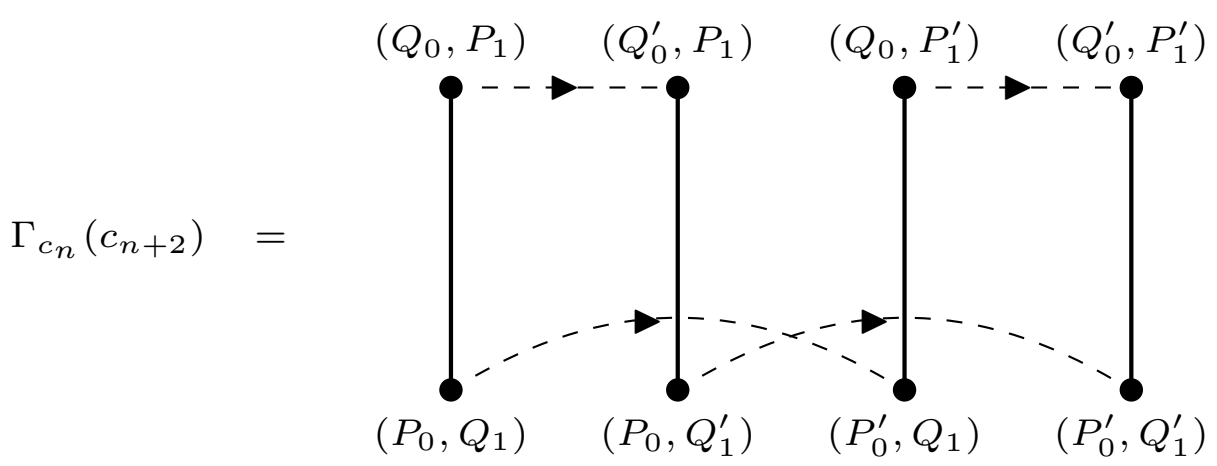
    \begin{center}
        \begin{tikzpicture}[font=\tiny]
            \node (a) at (-3,0) {$\Gamma_{c_n}(c_{n+2}) \;\;\; =$};
            \node (B) at (-1,-1) {};
            \filldraw[black] (-1,-1) circle (2pt);
            \node (b) at (-1,-1.3) {$(P_0,Q_1)$};
            \node (C) at (-1,1) {};
            \filldraw[black] (-1,1) circle (2pt);
            \node (c) at (-1,1.3) {$(Q_0, P_1)$} edge[-, thick] (b);
            \node (B2) at (0.5,-1) {} ;
            \filldraw[black] (0.5,-1) circle (2pt);
            \node (b2) at (0.5,-1.3) {$(P_0,Q_1')$};
            \node (C2) at (0.5,1) {} edge[-, dashed] (C) ;
            \node (C2C) at (-0.25,1) {$\blacktriangleright$};
            \filldraw[black] (0.5,1) circle (2pt);
            \node (c2) at (0.5,1.3) {$(Q_0', P_1)$} edge[-, thick] (b2);
            \node (B3) at (2,-1) {} edge[-, dashed, bend right] (B);
            \node (B3B) at (0.3,-0.55) {$\blacktriangleright$};
            \filldraw[black] (2,-1) circle (2pt);
            \node (b3) at (2,-1.3) {$(P_0',Q_1)$};
            \node (C3) at (2,1) {} ;
            \filldraw[black] (2,1) circle (2pt);
            \node (c3) at (2,1.3) {$(Q_0, P_1')$} edge[-, thick] (b3);
            \node (B4) at (3.5,-1) {} edge[-, dashed, bend right](B2);
            \node (B3B) at (1.8,-0.55) {$\blacktriangleright$};
            \filldraw[black] (3.5,-1) circle (2pt);
            \node (b4) at (3.5,-1.3) {$(P_0',Q_1')$};
            \node (C4) at (3.5,1) {} edge[-, dashed] (C3);
            \node (C4C3) at (2.75,1) {$\blacktriangleright$};
            \filldraw[black] (3.5,1) circle (2pt);
            \node (c4) at (3.5,1.3) {$(Q_0', P_1')$} edge[-, thick] (b4);
        \end{tikzpicture}
    \end{center}
    As we can see, we are left with a single component, and orienting these edges in the direction specified by the ordered tuples in $\M_{c_n}$, we see then that $F(C) = f(I_1) + f(I_2) + f(I_3) + f(I_4) =0$ and  $D(C)=2\equiv 0$. Hence $\sq^2_{c_n} (c_{n+2}) = 1$.
\end{example}
When applied to topological spaces, Steenrod squares are natural. In the case of framed 1-flow categories, we have the following result. 
\begin{lemma}[\cite{schuetz2022secondsteenrodsquaresodd}]
\label{lemma:sesConnHom}
    Let %
    \[
    0\longrightarrow \D \longrightarrow \Mm \longrightarrow \C\longrightarrow 0
    \]
    be a short exact sequence of framed 1-flow categories, and $\delta$ be the connecting homomorphism on homology. Then %
    \[
    \delta^* \circ \Sq^2 = \Sq^2 \circ \delta ^* : H^{k-1} (\C; \Z/2\Z) \to H^{k+2}(\D;\Z/2\Z).
    \]
\end{lemma}

We make some remarks about these special graph structures, which we include as a pair of lemmas. These will be useful later when showing that the Cartan formula for the product holds. 

\begin{lemma}\label{evenEdge}
    We have an even number of directed edges for a special graph structure with $\partial V = \emptyset$
\end{lemma}

\begin{proof}
    This follows from the fact that each component must be a cycle, and since each vertex is attached to an edge in $E'$, and an edge in $E-E'$, it must be an even cycle. First note that for any \textit{undirected edge}, which are all edges in $E-E''$, we have a sign change, ie for $e=\{v_1,v_2\}$, we have that $s(v_1) \neq s(v_2)$. There must be an even number of such edges, hence there is an even number of edges in $E''$ which are directed.
\end{proof}

A consequence of this is that if there are an odd number of intervals in some $\Gamma_\varphi(c^{k+2})$, then there will be an odd number of undirected edges in $E-E'$. 
\begin{lemma}\label{coboundaryFix}
     Let $\varphi$ be a cocycle in degree $k$ and $\{c_{h}^k\}$ , and let $\{c_i^{k+2}\}$ be the generators in $C^{k+2} (\C;\Ztwo)$ for which the SGS $\Gamma_\varphi (c_i^{k+2})$ contains an odd number of framed edges. Then $\sum c_i^{k+2}$ is a coboundary.
\end{lemma}
\begin{proof}
    We shall construct a cochain $\psi$ in degree $k+1$ with the property that generators $c^{k+2}$ which are in the image of $\delta (\psi)$ are exactly those whose SGS $\Gamma_\varphi (c_i^{k+2})$ have an odd number of framed edges. 

    Given some SGS for a cocycle representative $\varphi$, we have a combinatorial boundary matching $M=\{(P_i,P_i')\}_{i=1}^n$. We define $\psi$ to be the cochain in degree $k+1$ which counts the pairs $(P_i,P_i')$ with $s(P_i)\neq s(P_i')$ with multiplicity modulo 2. That is, for $P_i,P_i' \in \M(c_j,\varphi)$:
    \[
    \psi := \sum_{\substack{(P_i,P_i')\in M,\\ s(P_i)\neq s(P_i')}} c_j^{k+1}.
    \]
    See Figure \ref{fig:cmb1flow} for an example of one such pair, and note that for each pair $(P_i,P_i')$ there is a collection $\{Q_{i_k}\}_k = \M(c_i^{k+2},c_j^{k+1})$. When we have 
\[|\M(c_i^{k+2},c_j^{k+1})|\equiv 1 \;\mathrm{mod}\;2,
\]
the pair $(P_i, P_i')$ produces an odd number of undirected edges, hence, $\delta (\psi)$ will only contain $c_i^{k+2}$ in its sum if there are an odd number of undirected edges over all, since $\psi$ only contains generators with an odd number of the required pairs.
    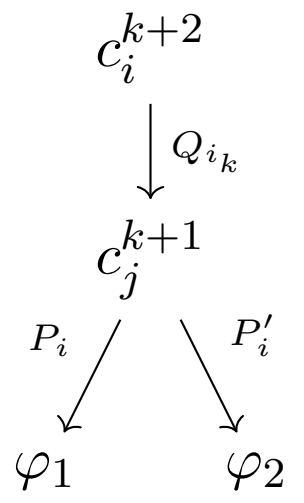
\begin{figure}[htbp]
        \centering
        \begin{tikzpicture}[every edge quotes/.style = {auto, font=\ttfamily\tiny,
      text=black, fill=none},scale = 0.7]
            \node (a) at (-1,0) {$\varphi_{1}$};
            \node(b) at (1,0) {$\varphi_{2}$}; 
            \node(c) at (0,2) {$c_j^{k+1}$} edge["$P_i$",swap, ->] (a) 
                                        edge["$P_i'$", ->](b);
            \node (d) at (0,4) {$c_i^{k+2}$} edge["$Q_{i_k}$", ->] (c);
        \end{tikzpicture}
        \caption{An example pair $(P_i,P_i')$ and 0-dimensional moduli spaces up to $c^{k+2}$}
        \label{fig:cmb1flow}
    \end{figure}
    Equivalently, suppose we have an odd number of undirected edges in $\Gamma_\varphi(c_j^{k+2})$. Then we have an odd number of pairs $((P,Q), (P',Q))$ with $(P,P')\in M_\varphi$ and $s(P)\neq s(P')$. For each such $Q$, we may associate an object in degree $k+1$ by $\target(Q)$, and taking the sum of these modulo $2$ will give a cochain associated to each special graph structure.
    
    Repeating this across all objects in degree $k+2$, we can take the \textit{union} of these cochains to form a $(k+1)$-cochain $\eta$. By construction then $\delta(\eta)$ will evaluate non-trivially exactly when there is an odd number of undirected edges in $\Gamma_\varphi(c_i^{k+2})$.
\end{proof}%
\section{Products and suspensions}
\label{sec:products}
\subsection{Products}

In this section we define a framed 1-flow category $\C\otimes \D$ and we justify that this is a suitable definition for a product.

\begin{definition}
\label{prodDef}
Let $(\C,s_1,f_1)$ and $(\D,s_2,f_2)$ be framed 1-flow categories. Then, for any $c\in \Ob(\C)$ and $d\in\Ob(\D)$ we define $(\C\otimes \D, s, f)$ so that the objects are given by 
\[\Ob{\Cprod} = \Ob{\C} \times \Ob{\D},\]
and moreover $|(c,d)| = |c|+ |d|$. The 0-dimensional moduli spaces are defined by
\[\M((c',d),(c,d))= \M_{\C}(c',c)\times \{d\} \text{ and,}\]
\[\M((c,d'),(c,d)) = \{c\}\times \M(d',d).\] 
We use the shorthand $(P,d) = Pd)$ and $(c,Q) = cQ$ and set $s(Pd) = s(P)$, $s(cQ)=|c| + s(Q)$. The 1-dimensional moduli spaces are defined by
        \begin{align}
            \morC{(c'',d)}{(c,d)} &= \morC[\C]{c''}{c} \times \{d\},\\
           \morC{(c',d')}{(c,d)} &= \morC[\C]{c'}{c}\times [0,1] \times \morC[\D]{d'}{d},\\
            \morC{(c,d'')}{(c,d)} &= \{c\} \times \morC[\D]{d''}{d}.
        \end{align}
We use the shorthand $P\otimes Q$ for intervals of the type (3.2). These are framed so that for circles $S\times \{d_l\}$, $\{c_k\}\times T$, we have,
            \begin{align*}
                f(S\times \{d\}) &= f_1(S)\\
                f( \{c\}\times T) &= f_2(T),
            \end{align*}
and for intervals $(I,d), P\otimes Q, (c,J)$ we have,
            \begin{align*}
                f((I,d)) &= f_1(I),\\
               f(P\otimes Q) &= s(P)(\gr_s(P) + s(Q)),\\
                f((c,J)) &= f_2(J) + (s(Q_1) + s(Q_2))|c|,
            \end{align*}
            where $\partial J = (Q_1,Q_3) \sqcup (Q_2,Q_4)$. 
\end{definition}

\begin{theorem}
    Let $(\C,s_1,f_1)$ and $(\D,s_2,f_2)$ be framed 1-flow categories. Then $(\Cprod , s, f)$ is a framed 1-flow category, and 
    \[
    C^n(\Cprod,\Z) \cong \bigoplus_{k+l=n} C^n(\C_1;\Z) \otimes C^l(\C_2;\Z).
    \]
\end{theorem}
\begin{proof}
    The final statement follows by construction. The choice of sign-assignment on the product yields the usual sign conventions on the tensor product of chain complexes so we have an identification $C_*(\C\otimes \D;\Z) = C_*(\C;\Z)\otimes C_*(\D;\Z)$.

    Since the sign assignment condition is automatically satisfied by the above, we must only check the framing of the 1-dimensional moduli spaces against their compatibility condition.
Let the subscripts on the objects indicate degrees in the respective flow categories. Then we have the following cases to consider:
    \begin{enumerate}[(i)]
        \item $\partial\morC{(c_{k+3},d_l)}{(c_{k},d_l)}$
        \item $\partial \morC{(c_{k+2},d_{l+1})}{(c_k,d_l)}$
        \item $\partial \morC{(c_{k+1},d_{l+2})}{(c_k,d_l)}$
        \item $\partial \morC{(c_k,d_{l_3})}{(c_k,d_l)}$
    \end{enumerate}
In case (i), it is clear that the framing values will be identical to those appearing in $\C$, and since the compatibility condition is satisfied in $\C$, it will be here as well. 

In case (ii), we may have intervals $I$ and circles, $C$, from $\morC[\C]{c_{k+2}}{c_k}$ and points $Q$ in $\morC[\D]{d_{l+1}}{d_l}$. For each pair of a point $Q$ and circle $C$ we have a pair of circle components $(C,d_{l+1})\times \{c_kQ\}$ and $\{c_{k+2}Q\}\times (C,d_0)$ in the boundary. Since the framing does not change, these cancel out modulo 2. For each interval and point pair, in the product we get a 3-cube where $\partial \morC{(c_{k+2},d_{l+1})}{(c_k,d_l)}$ a hexagon. The objects and morphisms are shown in general in Figure \ref{fig:case2cube}, and the moduli space is shown in Figure \ref{fig:case2hex}. Note that we may indeed have $c_{k+1} = c_{k+1}'$, and therefore either $P_1 = P_2$, or $P_3 = P_4$. In these cases, the moduli spaces look the same and will simplify in much the same way as our more general case in Figure \ref{fig:case2cube}.
% 
%Cube from I x Q
\begin{figure}[htbp]
    \centering
    \begin{tikzpicture}[every edge quotes/.style = {auto, font=\ttfamily\small,
      text=black, fill=none}, scale =0.8, transform shape]
% 
    %S
        \node (c2d1) at (3,6) {$(c_2,d_1)$};
        \node (c2d0) at (0,4) {$(c_2,d_0)$} edge ["$c_2 Q$",<-,sloped] (c2d1);
        \node (c1d1) at (3,4) {$(c_1,d_1)$} edge ["$ P_2d_1$",<-] (c2d1);
        \node (c11d1) at (6,4) {$(c_1',d_1)$}edge [ "$P_2'd_1$",<-,sloped] (c2d1);
        \node (c1d0) at (0,2) {$(c_1,d_0)$} edge ["$P_2d_0$", <-] (c2d0)
                                            edge [ "$c_1Q$",opacity =0.5,<-,sloped,near end] (c1d1);
        \node (c11d0) at (3,2) {$(c_1',d_0)$} edge [swap, "$P_2'd_0$", <-,sloped,near end] (c2d0)
                                              edge [swap,"$c_1'Q$", <-,sloped,near end] (c11d1);
        \node(c0d1) at (6,2) {$(c_0,d_1)$} edge ["$P_1d_1$",opacity =0.5,<-,sloped,near end] (c1d1)
                                            edge [swap, "$P_1'd_1$", <-] (c11d1);
        \node (c0d0) at (3,0) {$(c_0,d_0)$} edge [swap, "$P_1d_0$",<-,sloped] (c1d0)
                                            edge ["$P_1'd_0$", <-] (c11d0)
                                            edge [swap, "$c_0Q$", <-,sloped] (c0d1);
        %%%%%
\end{tikzpicture}
        \caption{The product 1-flow category arising from a `square' in $\C$ with four distinct objects and morphisms $P_1,P_2, P_3, P_4$, and a 0-dim moduli space $Q$ in $\D$.}
    \label{fig:case2cube}
\end{figure}
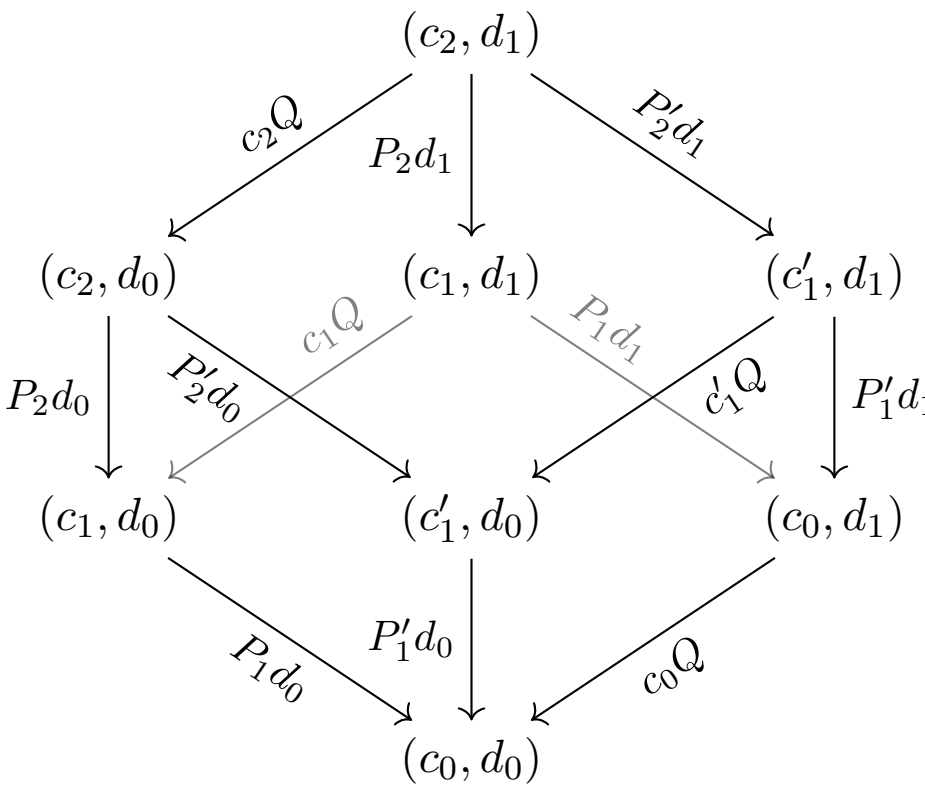
% 
%demonstrating the hexagon
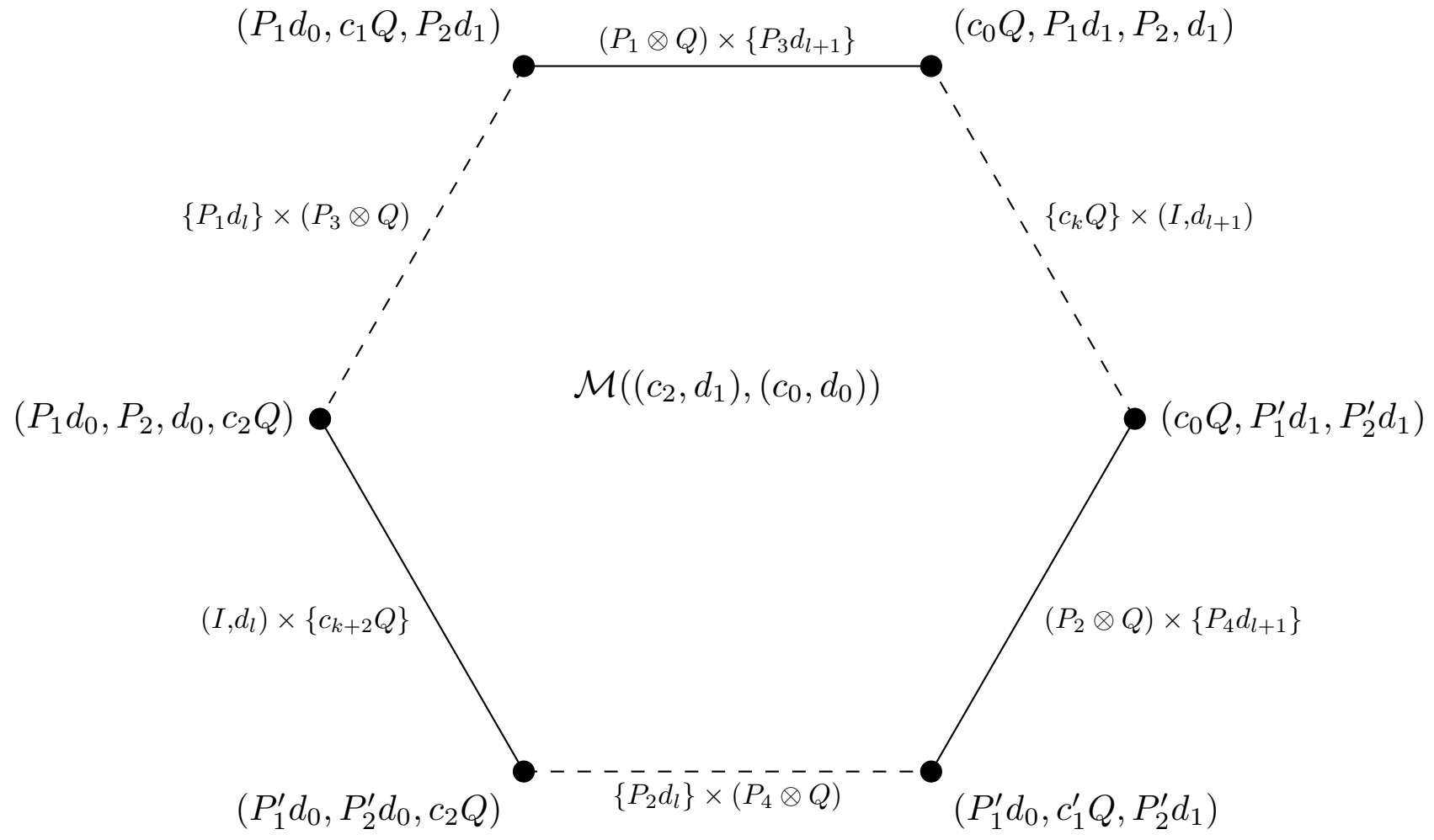
\begin{figure}[htbp]
    \centering
    \begin{tikzpicture}[scale=0.8, transform shape]
    \begin{scope}[font =\small]
        \newdimen\R
        \R=4cm
        \node[anchor = south] {$\M ((c_2,d_1),(c_0,d_0))$};
        \foreach \x/\l/\d in
            { 60/$\{c_kQ\} \times (I{,}d_{l+1})$/dashed,
            120/{$(P_1\otimes Q)\times \{P_3 d_{l+1}\}$}/,
            180/{$\{P_1d_l\} \times (P_3 \otimes Q)$}/dashed,
            240/{$(I{,}d_l) \times \{c_{k+2} Q\}$}/,
            300/{$\{P_2 d_l\} \times (P_4\otimes Q)$}/dashed,
            360/{$(P_2\otimes Q) \times \{P_4 d_{l+1}\}$}/%
            }
        \draw[postaction={decorate}] ({\x-60}:\R) edge[-,\d] node[auto,swap,scale = 0.8]{\l} (\x:\R);
        \foreach \x/\l/\p in
            { 60/{($c_0Q, P_1d_1, P_2,d_1$)}/above right,
            120/{($P_1d_0, c_1 Q, P_2 d_1$)}/above left,
            180/{$(P_1d_0, P_2,d_0,c_2Q)$}/left,
            240/{$(P_1'd_0, P_2'd_0, c_2Q)$}/below left,
            300/{$(P_1'd_0, c_1' Q, P_2'd_1)$}/below right,
            360/{$(c_0 Q, P_1'd_1,P_2'd_1)$}/right
            }
        \node[inner sep=2pt,circle,draw,fill,label={\p:\l}] at (\x:\R) {};
    \end{scope}
    \end{tikzpicture}
    \caption{The moduli space, captioned in centre, has a boundary hexagon which shows up in our compatibility condition.}
    \label{fig:case2hex}
\end{figure}
Reading framings counter-clockwise from $(c_0Q, P_1'd_1,P_2'd_1)$ in Figure \ref{fig:case2hex}, this hexagon contributes the following sum to the compatibility condition:
    \begin{align*}
    &\sum \tilde{f}(I\times\{Q\}) + \sum \tilde{f}(\{P\}\times J) =(1 + s(c_0 Q) + f((I,d_1))\\ &+ f(P_1\otimes Q) + (1 + s(P_1d_0) + f(P_2\otimes Q))+f((I,d_0))\\
    &  + (1 + s(P_1' d_0) + f(P_2'\otimes Q)) + f(P_1'\otimes Q).\\
    \end{align*}
By our definitions can be shown to sum to $1$, and therefore satisfies the compatibility condition.

In the next case we have a 0-dimensional moduli space in $\C$ and a 1-dimensional moduli space in $\D$. This works out in the same way and we have
\begin{align*}
&\sum \tilde{f}(I\times\{Q\}) + \sum \tilde{f}(\{P\}\times J) = (1+s(Pd_0)+ f((c_1,J)))\\&+ f(P\otimes Q_1) + (1 + s(c_0 Q_1) + f(P \otimes Q_2))+ f((c_0,J))\\
&+ (1 + s(c_0 Q_1') + f(P \otimes Q_2')) + f(P\otimes Q_1')\\
&= 1 + s(P) + s(P)(s(Q_1) + s(Q_2) + s(Q_1') + s(Q_2'))\\
&=1.
\end{align*}
Lastly, we have the case where the boundary of the 2-dimensional moduli space is directly from $\D$, which we depict in Figure \ref{fig:case3Prod}. Unlike the first case, our framing and sign assignment changes and the compatibility condition does not follow as obviously. Without any information on the moduli spaces, we must work with an arbitrary $2n$-gon. We will use $(A,B,C)\in \M(d_1,d_0)\times \M(d_2,d_1)\times \M(d_3,d_2)$.
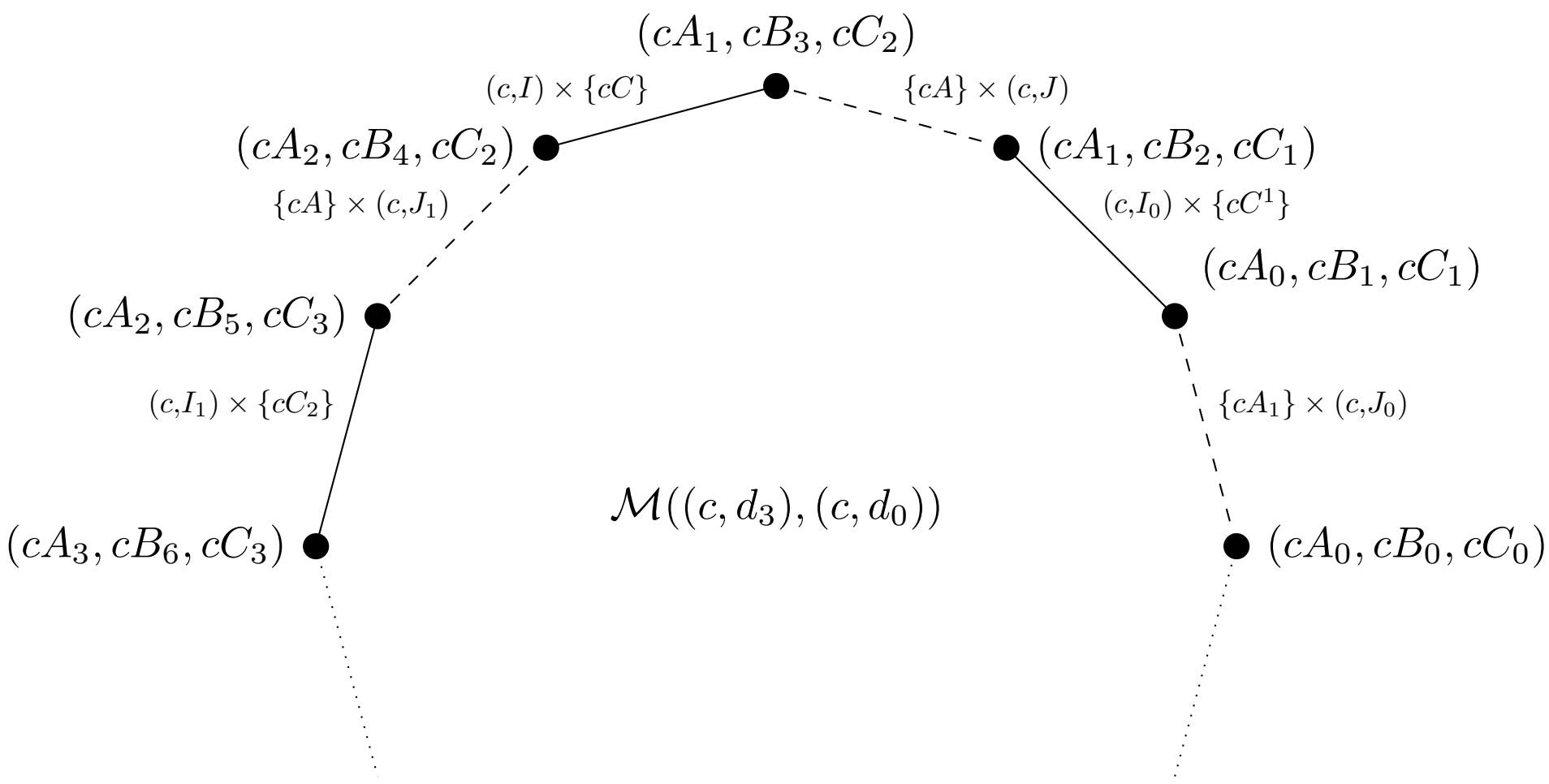
\begin{figure}[htbp]
    \centering
    \begin{tikzpicture}[font=\tiny, scale=0.75]

        \begin{scope}[font =\small]
            \newdimen\R
            \R=5cm
            \node[anchor = south] {$\M ((c,d_3),(c,d_0))$};
            \foreach \x/\l/\d in
                { 0//dotted,
                30/{$\{cA_1\}\times (c{,}J_0)$}/dashed,
                60/{$(c{,}I_0) \times\{cC^1\}$}/,
                90/{$\{cA\}\times (c{,}J)$}/dashed,
                120/{$(c{,}I)\times \{cC\}$}/,
                150/{$\{cA\}\times (c{,}J_1)$}/dashed,
                180/$(c{,}I_1)\times \{cC_2\}$/,
                210//dotted%
                }
            \draw[postaction={decorate}] ({\x-30}:\R) edge[-,\d] node[auto,swap,scale = 0.7]{\l} (\x:\R);
                \foreach \x/\l/\p in
                    {180/{$(cA_3,cB_6,cC_3)$}/left,
                    150/{$(cA_2,cB_5,cC_3)$}/left,
                    120/{$(cA_2,cB_4,cC_2)$}/left,
                    90/{$(cA_1,cB_3,cC_2)$}/above,
                    60/{$(cA_1,cB_2,cC_1)$}/right,
                    30/{($cA_0,cB_1,cC_1)$}/above right,
                    360/{$(cA_0,cB_0,cC_0)$}/right
                    }
                \node[inner sep=2pt,circle,draw,fill,label={\p:\l}] at (\x:\R) {};
        \end{scope}
\end{tikzpicture}
    
    \caption{Example boundary of  
    $\M ((c,d_3),(c,d_0))$.}
    \label{fig:case3Prod}
\end{figure}
Starting with the solid edges we have per edge \[\tilde{f}((c,I) \times \{cC\}) = (s(A_0)+s(A_1))|c| + f(J) = s(A_0)|c| + s(A_1)|c| + f(J).\] recalling the definition of $\Tilde{f}$ from Def. \ref{framingdef}. That is to say, in the final sum we can attribute to each vertex $(A,B,C)$  $+s(A)|c|$.

Since there are always two vertices of the form $(A, *, * )$ (those on either side of a dashed edge labelled $\{A\} \times J$) the framing on the solid edges contributes $f(I)$ to the sum, the same as they do in $\D$.

To deal with the dashed edges we make a simple reformulation and then reapply this process. For example, by definition $\tilde{f}(\{cA\}\times cJ) = 1 + s(cA) + (s(cB) + s(cB^1))|c| + f(J)$, but since the sign assignments sum to 1 across squares, we can change this to  $ 1 + s(A) + |c| + (s(C) + s(C^1) + 1)|c| + f(J) = 1 + s(A) + f(J) + s(C)|c| + s(C^1) |c|$, and repeating the above step for the solid edges we now see that since there are always two vertices of the form $(*,*,Z)$, we have the same sum in the compatibility condition as in $\D$.

This argument assumes that adjacent to a solid edge are two distinct dashed edges, and so it does not resolve the possibility of digons. These arise in one case, where $\partial \M(a,d) $ is equal to $I\times \{D\}\sqcup \{A\}\times J$, and hence
\[\M(d_{l+1}, d_l)\times \M(d_{l+2},d_{l+1})\times \M(d_{l+3},d_{l+2}) = (D_1,D_2,D_3)\sqcup (D_1,D_2',D_3)\]\
, with $s(D_2)\neq s(D_2')$. In this case, the sum is more straightforward, and we have
\begin{align*}
&f((c_k,I)) + 1 + s(c_k D_1) + f((c_k,J)) = (s(D_1)+s(D_1))|c_k| + f(I) \\ &+ 1 + s(D_1) + |c_k| + (s(D_2) + s(D_2'))|c_k| + f(J)\\
&= f(I) + 1 + s(D_1) + f(J),
\end{align*}
which is the same contribution to the compatibility condition from this component as in $\D$.

So, since the sign assignments and framings satisfy the desired equalities, $(\C\otimes \D, s,f)$ is framed.
\end{proof}

In order to justify our choices, we will show that the product behaves as expected with the second Steenrod square. 
\begin{theorem}[Cartan Formula] \label{thm:cartan}
    Let $(\C,s,f)$ and $(\D,s',f')$ be framed 1-flow categories, and let 
    $$i: H^k(\C; \Ztwo)\otimes H^l(\D ; \Ztwo) \rightarrow H^{k+l}(\C\otimes \D ; \Ztwo)$$
    be the inclusion from the Kunneth formula. Then we have:
    \[
    \Sq^2 (i(z_1 \otimes z_2)) = i(\Sq^2(z_1)\otimes z_2 + \Sq^1(z_1)\otimes \Sq^1(z_2) + z_1 \otimes \Sq^2(z_2))
    \]
    for all $z_1 \in H^k (\C; \Ztwo)$ and $z_2 \in H^l(\Ztwo)$.
\end{theorem}
\begin{proof}
    The statement concerns cohomology classes $z_1 \in H^k(\C)$ and $z_2\in H^l(\D)$ for which we can associate representative cochains $\varphi \in C^k(\C;\Z)$ and $\psi \in C^l(\D)$, which we further identify as a sum of generators given by the objects in $\C$ and $\D$. 

    Due to the isomorphism of the underlying cochain complexes of the product with the tensor product of cochain complexes, the cochain $\varphi\otimes \psi \in \C^k(\C;\Z)\otimes C^l(\D;\Z)$ can be considered as a cochain in $C^n (\C\otimes \D;\Z)$. Likewise, we can identify maps $\sq^2_\varphi\otimes \psi$, $\sq^1_\varphi \otimes \sq^1_\psi$, and $\varphi \otimes \sq^2_\psi$ via this isomorphism as cochains in $C^{n+2}(\C\otimes \D)$. We compare the output of $\sq^2_{\varphi\otimes \psi} \in C^{n+2}(\C\otimes \D)$ with the outputs of the previous three $(n+2)$-cochains on the following objects in $\C\otimes\D$, once again, labelling the grading in the subscripts:
    \begin{enumerate}[(1)]
        \item $(c_{k+2},d_l)$,
        \item $(c_{k+1},d_{l+1})$,
        \item $(c_k,d_{l+2})$, or
        \item otherwise, i.e $(c_m, d_n)$, where $m+n = k+l+2$, but not as above.
    \end{enumerate}
   Indeed, the cochain group $C^{n+2}(\C\otimes \D)$ splits as a direct sum $C_1 \oplus C_2 \oplus C_3 \oplus C_4$, where these are generated by objects of the forms in (1), (2), (3) and (4), respectively. From here it is clear that $\sq^2_\varphi\otimes \psi$ will only evaluate as non-zero on $C_1$, and likewise for $\sq^1_\varphi\otimes\sq^1_\psi$ on $C_2$ and $\varphi \otimes \sq^2_\psi$ on $C_3$. Provided these agree with $\sq^2_{\varphi\otimes \psi}$ on their corresponding direct summands, their sum will then necessarily agree across $C^{n+2}(\C\otimes \D)$. 
    
    Dealing with case (4) first, it is clear that the combinatorial Steenrod square will be trivially zero on these objects as we have $\M((c_m,d_n) , \varphi\otimes \psi) = \emptyset$.

    Since all 0-dimensional moduli spaces in $\C\otimes \D$ are of the form $P\times \{d\}$ or $\{c\}\times Q$, a combinatorial boundary matching on the 0-dimensional moduli spaces $\M_\C(c_{k+1},\varphi)$ and $\M_\D(d_{l+1}, \psi)$  naturally gives us a combinatorial boundary matching on the 0-dimensional moduli spaces given by $\M((c_{k+1},d_l),\varphi\otimes\psi)$ and $\M((c_k,d_{l+1}),\varphi\otimes\psi)$.

    To see how $\sq^2$ evaluates on objects of the form in (1), we note that the moduli space $\M((c_{k+2},d_l),\varphi\otimes\psi)$ is equal to $\M_\C (c_{k+2},\varphi)$, and with the same combinatorial matching we get identical special graph structures. The framing of intervals $I\times\{d_l\}$ and circles $S\times\{d_l\}$, and the sign assignments of objects $P\times \{d_l\}$ in this graph remain unchanged, hence $\sq^2_{\varphi\otimes\psi}(c_{k+2},d_l) = sq^2_\varphi(c_{k+2})\cdot\psi (d_l)$.

    To deal with (2), we must investigate the product between $\M_\C(c_{k+1},\varphi)$ and $\M_\D(d_{l+1},\psi)$, however, given the formula for $\sq^1$, which only counts the number of equivalently framed points in the combinatorial matching, we need only investigate the contribution of one of these pairs, and then from here we may simply sum over the contributions of each pair.  For each of these matched points, we have within our 1-flow category, the diagrams in Figure \ref{sq1pf} appearing, where we may have $c_k = c_k'$ or $d_l = d_l'$. 
\begin{figure}[htbp]
    \centering

    \begin{tikzpicture}[scale=1.5, font = \small]

        \begin{scope}
            \node (a) at (2,0) {$c_k$};
            \node (b) at (4,0) {$c_k'$};
            \node (c) at (3,1) {$c_{k+1}$} edge [swap, "$P$", ->] (a)
                                           edge [ "$P'$", ->] (b); 
        \end{scope}
        \begin{scope}[shift = {(4,0)}]
            \node (a) at (2,0) {$d_l$};
            \node (b) at (4,0) {$d_l'$};
            \node (c) at (3,1) {$d_{l+1}$} edge [swap, "$Q$", ->] (a)
                                           edge [ "$Q'$", ->] (b);
        \end{scope}
 
    \end{tikzpicture}
    
    \caption{These objects and morphisms appear in $\C$ and $\D$ for each pair of points in $\M(c_{k+1}, \varphi)$ and $\M(d_{l+1},\psi)$.}
    \label{sq1pf}
\end{figure}
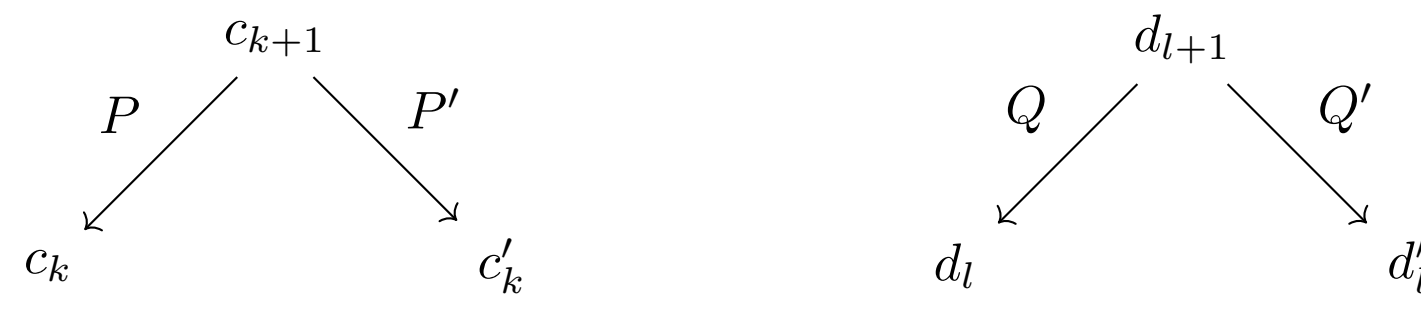
    We first consider a product of these moduli spaces, and then evaluate the Steenrod squares.

    Taking the product gives us a flow category as in Figure \ref{prodSq1} which can be seen as a full sub-category of $\C\otimes\D$.
    \begin{figure}[htbp]
        \centering

        \begin{tikzpicture}[font=\tiny, transform shape, scale=1.1]
        \node (a) at (0,3) {$(c_{k+1},d_{l+1})$};
        \node (b) at (-3,1.5) {$(c_{k+1},d_l)$} edge [sloped,thick, near end, "$c_{k+1} Q$",<-] (a);
        \node (c) at (-1,1.5) {$(c_{k},d_{l+1})$} edge [near start, sloped, thick,  "$Pd_{l+1}$", <-] (a);
        \node (d) at (1,1.5) {$(c_{k}',d_{l+1})$}  edge [sloped, thick, near start, "$P'd_{l+1}$",<-] (a);
        \node (e) at (3,1.5) {$(c_{k+1},d_{l}')$}  edge [sloped, thick, near end, "$c_{k+1}Q'$",<-] (a);
        \node (f) at (-3,0) {$(c_{k},d_{l})$} edge [thick, midway, "$Pd_l$", <-] (b)
                                                edge [ sloped,thick, near end, "$c_kQ$", <-] (c);
        \node (g) at (-1,0) {$(c_{k}',d_{l})$} edge [sloped, opacity = 0.6,near start, "$P'd_l$", <-] (b)
                                                edge [sloped,opacity = 0.6, very near start, "$c_k'Q$", <-] (d);
        \node (h) at (1,0) {$(c_{k},d_{l}')$} edge [sloped, thick, very near start,  "$Pd_l'$", <-] (e)
                                                edge [sloped, thick, near start, "$c_kQ'$",  <-] (c);
        \node (i) at (3,0) {$(c_{k}',d_{l}')$} edge [swap, thick,"$P'd_l'$",  <-] (e)
                                                edge [sloped, opacity =0.6, "$c_k'Q'$", near start, <-] (d);
            
        \end{tikzpicture}
        
        \caption{The product of the 1-flow categories in Figure \ref{sq1pf}.}
        \label{prodSq1}
    \end{figure}
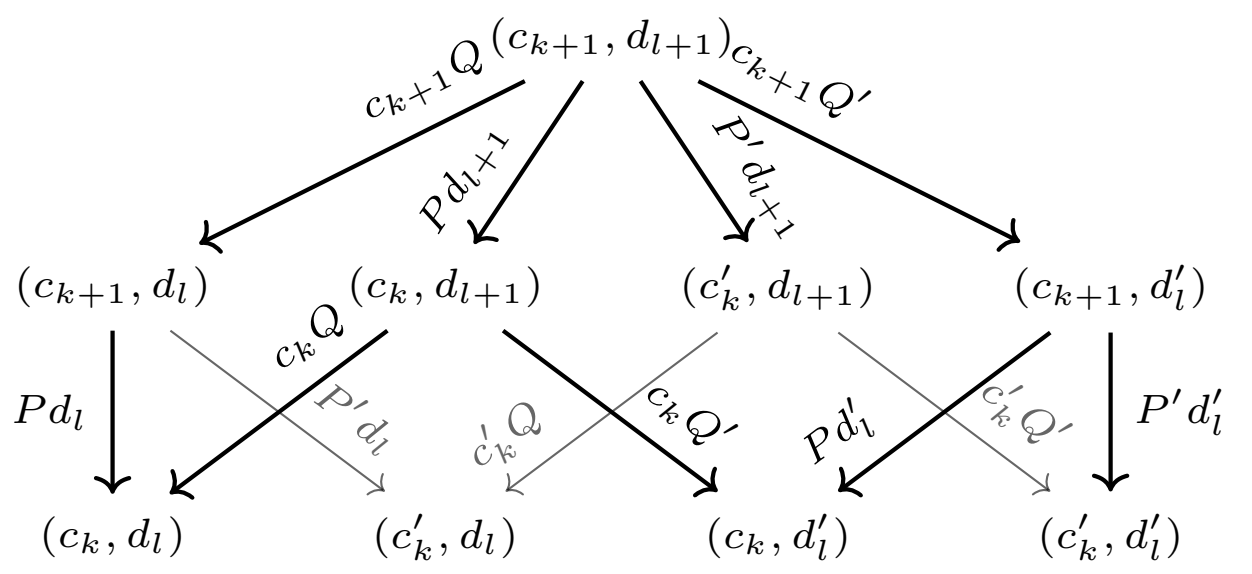
The moduli spaces $\M((c_{k+1},d_{l+1}), \varphi\otimes \psi)$ encoded in our special graph structure will have four intervals $I_1,I_2,I_3, I_4$, which correspond respectively to the moduli spaces from $(c_{k+1},d_{l+1})$ to $(c_k,d_l), (c_k,d_l'), (c_k'd_l'), (c_k',d_l)$. The special graph structure will take the form in Figure \ref{sgs1}.

\begin{figure}[htbp]
    \centering

    \begin{tikzpicture}[font = \tiny, scale=0.85]
        \node[anchor = east] (a) at (-1,4) {$(Pd_l, c_{k+1}Q)$};
        \node[anchor = west] (b) at (1,4) {$(c_kQ, Pd_{l+1})$} edge [thick, midway, "$s(P)(k+1+s(Q))$"] (a);
        \node[anchor = east]  (c) at (-1,2) {$(Pd_l', c_{k+1}Q')$};
        \node[anchor = west]  (d) at (1,2) {$(c_kQ', Pd_{l+1})$} edge [thick, midway, "$s(P)(k+1+s(Q'))$"] (c) edge [dashed, bend right=90] (b);
        \node[anchor = east]  (e) at (-1,0) {$(P'd_l', c_{k+1}Q')$}edge [dashed] (c);
        \node[anchor = west]  (f) at (1,0) {$(c_k'Q', P'd_{l+1})$} edge [thick, midway, "$s(P')(k+1+s(Q'))$"] (e) ;
        \node[anchor = east]  (g) at (-1,-2) {$(P'd_l, c_{k+1}Q)$}edge [dashed,opacity = 0.6, bend left=90] (a);
        \node[anchor = west]  (h) at (1,-2) {$(c_k'Q, P'd_{l+1})$} edge [thick, midway, "$s(P')(k+1+s(Q))$"] (g) edge [dashed, bend right=90] (f) ;
        
    \end{tikzpicture}
    
    \caption{The intervals $I_i$ in $\M((c_{k+1},d_{l+1}),\varphi\otimes\psi)$ with annotated framing values on the edges}
    \label{sgs1}
\end{figure}
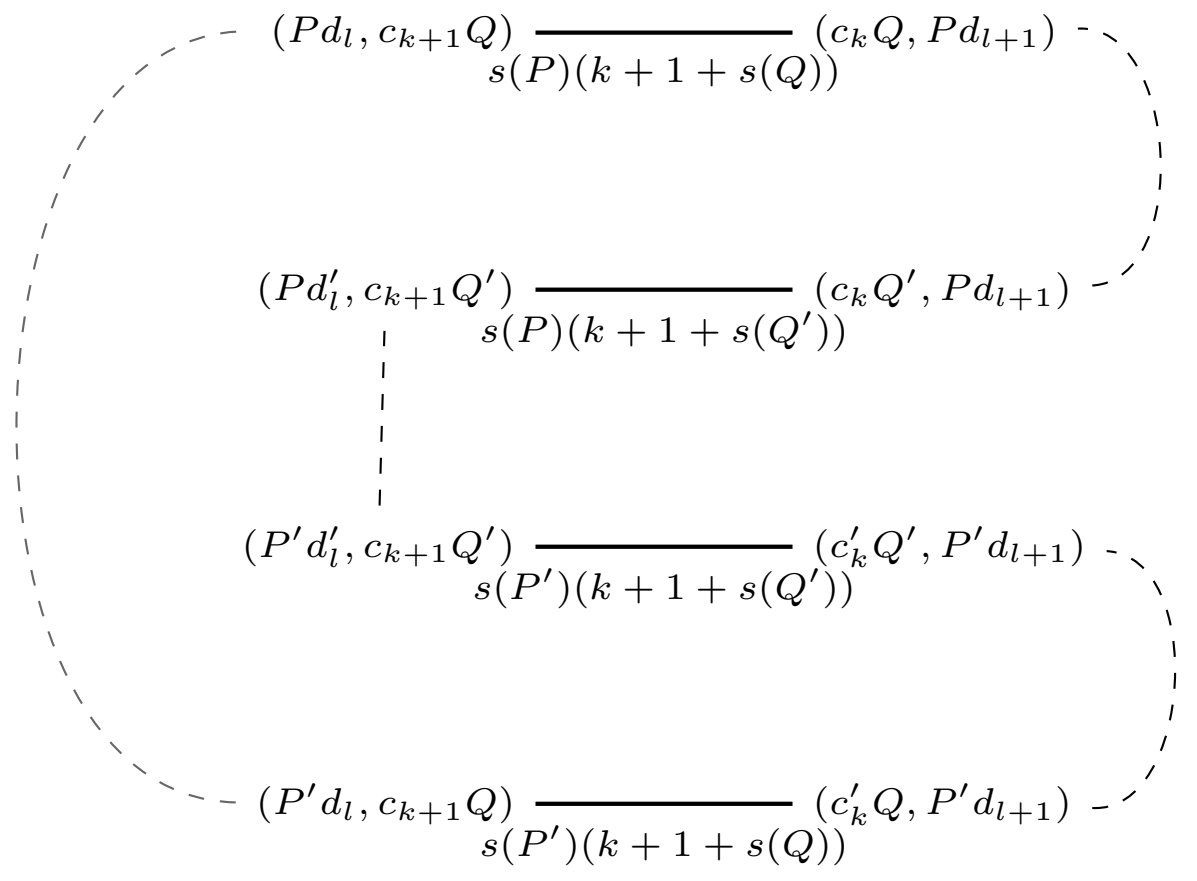
Note, if $s(P) = s(P')$, then the dashed edges corresponding to that matched pair will be marked with arrows, and we contribute $+1$ to the sum, and likewise if $s(Q) = s(Q')$. We may then contribute $(s(P)+s(P') + 1)$ and $(s(Q) + s(Q') +1)$ to the sum to account for this. Hence, we have:
\begin{align*}
&sq^2_{\varphi\otimes\psi}((c_{k+1},d_{l+1})) = 1 + (s(P)+s(P')+1) + (s(Q)+s(Q')+1) +\\
&s(P)(1+k+s(Q)) + s(P')(1+k+s(Q')) + s(P')(1+k+s(Q') +\\
&s(P')(1+k+s(Q)).
\end{align*}
Which simplifies to 
$$sq^2_{\varphi\otimes\psi}((c_{k+1},d_{l+1})) = (s(P)+s(P')+1)(s(Q)+s(Q')+1).$$
This component of the special graph structure evaluates to $1$ if and only if both $s(P) = s(P')$ and $s(Q) = s(Q')$. Summing these contributions across all matched pairs is equivalent to evaluating the combinatorial $sq^2$ on each component of $\Gamma_\varphi$ and summing these contributions up. On the other hand, evaluating $\sq^1_\varphi\otimes\sq^1_\psi$ gives the same outcome; the product of the sums corresponding to the combinatorial first Steenrod square will be the sum of the occasions when both $s(P)=s(P')$ and $s(Q)=s(Q')$, hence,
$$sq^2_{\varphi\otimes\psi}((c_{k+1},d_{l+1})) = sq^1_\varphi(c_{k+1})\otimes sq^1_\psi(d_{l+1}).$$
Lastly, evaluating $sq^2_{\varphi\otimes \psi}((c_k,d_{l+2}))$, we find some discrepancies between the left and right hand side of the equation. We examine a part of the special graph structure corresponding to $sq^2_{\varphi\otimes\psi}((c,d_{l+2}))$ in Figure \ref{sq2iii} to demonstrate the difference. Compared with the graph structure corresponding to $sq^2_\psi(d_{l+2})$, the difference is only seen in the framing of the intervals, where we have an additional summand of $s(Q_1)|c| + s(Q_1')|c|$. The directed edges connected to these framed edges have $s(Q_1) = s(Q_2)$, hence these will cancel, but for each undirected edge we end up with an extra $|c|$. Hence, when there are an odd number of undirected edges, we have an inequality. However, by Lemma \ref{coboundaryFix} we know that the two representatives must then differ by a coboundary, hence the left and right hand sides agree up to cohomology.
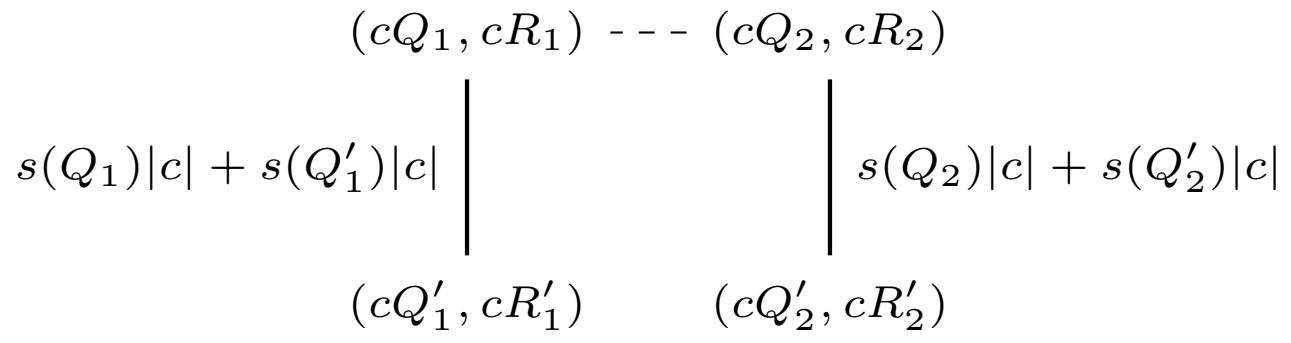
\begin{figure}[htbp]
    \centering
    
    \begin{tikzpicture}[font = \tiny, scale =1.2, transform shape]
        \node (a) at (-1,0.75) {$(cQ_1, cR_1)$};
        \node(b) at (-1,-0.75) {$(cQ_1', cR_1')$} edge [thick, midway, "$s(Q_1)|c| + s(Q_1')|c|$"] (a);
        \node (c) at (1,0.75) {$(cQ_2, cR_2)$} edge [dashed] (a);
        \node (d) at (1,-0.75) {$(cQ_2', cR_2')$} edge [thick, midway, "$s(Q_2)|c| + s(Q_2')|c|$",swap] (c);
    \end{tikzpicture}
    \caption{An example of a section of a special graph structure corresponding to $sq^2(\varphi\otimes\psi)[(c_k,d_{l+2})]$ with labelling on edges corresponding to the difference in framing of this interval in $\D$ and $\C \otimes\D$.}
    \label{sq2iii}
\end{figure}
\end{proof}
In the homotopy category of spectra, the smash product is, up to equivalence, a commutative, associative product with unit given by the sphere spectrum $\mathbb{S}$ \cite{Adams_1995}. Following from the Cartan formula, further justification of our product can be seen from the following lemma.

\begin{lemma}
\label{lemma:assocProd}
    The following properties hold:
    \begin{enumerate}[(i)]
        \item There exist a 1-flow category $\Sp_0$ which acts as the identity with respect to the product.
        \item The product is associative.
    \end{enumerate}
\end{lemma}

\begin{proof}
We define $\Sp_0$ to simply be the framed 1-flow category with a unique object in degree 0, and it is clear that this acts as an identity element, as the framed 1-flow categories $\C$, $\Sp_0\otimes \C$ and $\C\otimes \Sp_0$ are equal on the levels of objects, moduli spaces and framings. 

We omit the proof of associativity, as this only requires careful application of the definitions.
\end{proof}%
Due to the asymmetric way that intervals and points are framed in the product- the product will not commute directly. Instead, it might be desirable to understand framed 1-flow categories up to some notion of homotopy equivalence, where it should be expected that the product commutes up to homotopy (due to the Theorem \ref{thm:cartan}, it can be shown that the Steenrod square data is the same, for instance). We do not discuss in this paper any notions of equivalence for framed 1-flow categories, and leave this to be explored in the author's thesis, \cite{thesis}.  
\section{Mapping cones and the second Steenrod square}
We begin by defining the suspension of a framed 1-flow category.
\begin{definition}
    Let $\C$ be a 1-flow category, and $\Sp_1$ the 1-flow category with a single object in degree 1. We define the suspension of $\C$ by $\Sigma \C=\Sp_1\otimes \C$.
\end{definition}
Notably this is defined differently in \cite{schuetz2022secondsteenrodsquaresodd}, where Sch\"utz's suspension may be identified with $\C\otimes \Sp_1$. Intervals with boundary $(P,Q)\sqcup (P',Q')$ change their framing by an additional $(s(P) + s(P'))$, and the sign assignment across the 1-flow category becomes $s_{\Sigma\C}(P)= s(P) + 1$.

\begin{definition}
    Given framed 1-flow categories $\C$, $\D$, $\Mm$, we say $\Mm$ is a mapping cone framed 1-flow category if we have $\Ob(\Mm) = \Ob(\Sigma \C) \sqcup \Ob(\D)$, and moreover, $\Mm$ fits into a short exact sequence of framed 1-flow categories, with $\Sigma \C$ and $\D$ fitting into $\Mm$ as subcategories in the obvious way:
     \[
        0\longrightarrow \D \longrightarrow \Mm \longrightarrow \Sigma \C \longrightarrow 0
    \]
    We may also write $\Mm_f$, where $f:C_*(\C;\Z)\to C_*(\D;\Z)$ is the chain map one gets by considering the differential of the underlying chain complex of the mapping cone.
\end{definition}
Using mapping cones, we provide an alternative definition of the combinatorial second Steenrod square, leveraging the naturality statement of Lemma \ref{lemma:sesConnHom} to ensure that once the appropriate compatibility conditions are satisfied we have a well-defined definition.

\begin{definition}
    Let $\C$ be a framed 1-flow category. We define the $n$-truncated framed 1-flow category $\C^{(n)}$ to be the full subcategory of $\C$ which includes all objects $x$ of $\C$ with $|x|\leq n$.
\end{definition}

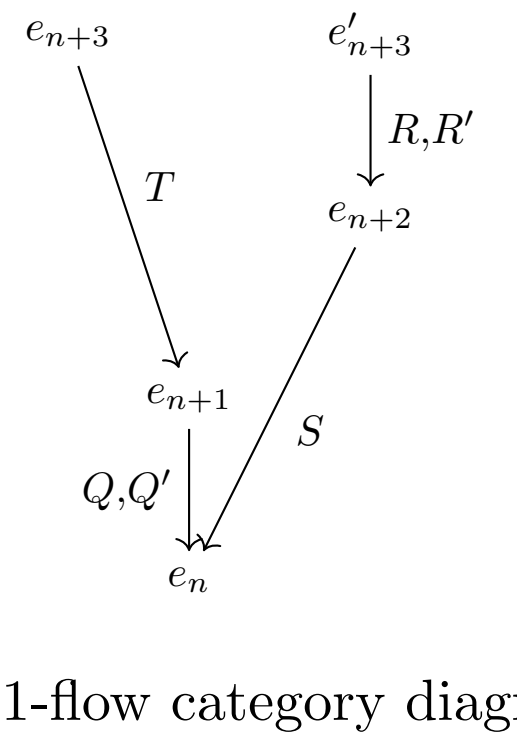
\begin{figure}[htbp]
    \centering
    \begin{tikzpicture}[font = \small, scale=0.8, transform shape]
        \node (a) at (0,0) {$e_{n}$};
        \node (b) at (0,1.5) {$e_{n+1}$} edge [->, swap, "$Q{,}Q'$"] (a);
        \node (c) at (1.5,3) {$e_{n+2}$} edge [->, "$S$"] (a);
        \node (d) at (-1,4.5) {$e_{n+3}$} edge [->, "$T$"] (b);
        \node (e) at (1.5,4.5) {$e_{n+3}'$} edge [->, "$R{,}R'$"] (c);       
    \end{tikzpicture}
    \caption{1-flow category diagram for $\mathcal{K}_n$}
    \label{fig:k2n}
\end{figure}

\begin{definition} \label{def:eilenberg}
    Let $(\mathcal{K_n},s,f)$ be the framed 1-flow category with $\Ob(\mathcal{K}_n) = \{e_n, e_{n+1}, e_{n+2}, e_{n+3}, e_{n+3}'\}$, (where the gradings are listed in the subscripts). Let $S,T$ be framed circles so that $f(S)=f(T)=0$, and $Q,Q',R,R'$ be framed points so that $s(Q)=s(Q')=s(R)=s(R') =0$. Then let the moduli spaces of $\mathcal{K}_n$ be given by $\M(e_{n+1},e_n) = \{Q,Q'\}$, $\M(e_{n+2} ,e_n) = S$, $\M(e_{n+3},e_{n+1}) = T$, and $\M(e_{n+3}', e_{n+2}) = \{R,R'\}$, with all other possible moduli spaces empty.
    
    A diagram for $\mathcal{K}_n$ is shown in \ref{fig:k2n}.
\end{definition}

\begin{remark} 
We think of this framed 1-flow category as a framed 1-flow category for the $k+3$ truncation of the stable Eilenberg-Maclane space $\oldcal{K}(\Z/2\Z,n)$. 
\end{remark}

\begin{proposition}\label{prop:altDefSq2}
    Let $\C$ be a framed 1-flow category, $\varphi\in C^n(\C;\Ztwo)$ be a cocycle represented by a collection of objects $\{x_i\}$ in degree $n$. There exists a framed mapping cone 1-flow category $(\Mm_\varphi,s,f)$ which lies in a short exact sequence of framed 1-flow categories
    \[
    0\longrightarrow \mathcal{K}_n \longrightarrow \Mm_\varphi \longrightarrow \Sigma (\C^{(n+3)}) \longrightarrow 0
    \] 
    Furthermore, if we express the differential on $\Mm_\varphi$ acting on the underlying chain complex as
    \[
    d_{\Mm_f} : =  \begin{bmatrix}
                        -d_\C &0\\
                        u &d_{\mathcal{K}_n}
                \end{bmatrix}   ,
    \]
    so that $u$ is a chain map $u: C_*(\C;\Z)\to C_*(\mathcal{K}_n;\Z) $ then we have on cohomology $u^*([e_n])=[\varphi]$, $u^*([e_{n+1}]) = \Sq^1 ([\varphi])$ and $u^*([e_{n+2}]) = \Sq^2([\varphi])$.
\end{proposition}

\begin{proof}
    In order to show that this is framed, we need to define the moduli spaces 
    We start by picking for $\varphi \in C^k(\C;\Z/2\Z)$ a combinatorial matching $M$ of the moduli spaces $\M(c_j, \varphi)$. This pairing will become essential for deciding where the intervals go in the moduli spaces $\M(\Sigma c_j, e_n)$.

    The objects of $\Mm_\varphi$ are defined to be simply a union of those in $\Sigma(\C^{(n+3)})$ and those in $\mathcal{K}_n$, and as implied by the short exact sequence of framed 1-flow categories, the degrees, the moduli spaces, and the framings internal to $\Sigma(\C^{(n+3)}), \mathcal{K}_n$ are preserved. 
    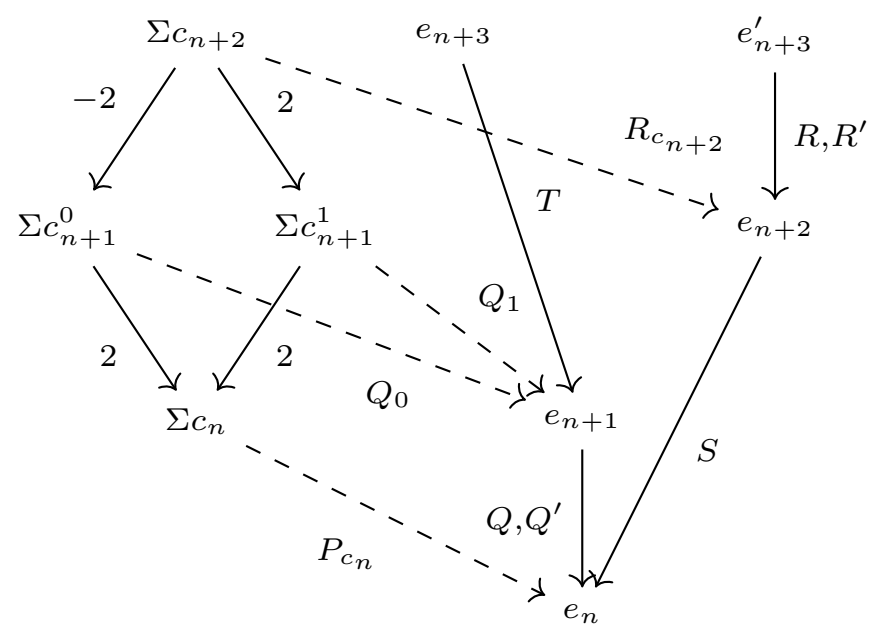
\begin{figure}
        \centering
            \begin{tikzpicture}[font=\tiny, scale =0.85, transform shape]
                \node (a) at (0,0) {$e_{n}$};
                \node (b) at (0,1.5) {$e_{n+1}$} edge [->, swap, "$Q{,}Q'$"] (a);
                \node (c) at (1.5,3) {$e_{n+2}$} edge [->, "$S$"] (a);
                \node (d) at (-1,4.5) {$e_{n+3}$} edge [->, "$T$"] (b);
                \node (e) at (1.5,4.5) {$e_{n+3}'$} edge [->, "$R{,}R'$"] (c); 
                \node (u) at (-3,1.5) {$\Sigma c_n$} edge [->, dashed, swap,"$P_{c_n}$"] (a);
                \node (v) at (-2,3) {$\Sigma c_{n+1}^1$} edge [->, "$2$"] (u)
                                                        edge [->, dashed, "$Q_1$"] (b);
                \node (w) at (-4,3) {$\Sigma c_{n+1}^0$} edge [->, "$2$", swap] (u)
                                                         edge [->, dashed, "$Q_0$",swap, near end] (b);
                \node (x) at (-3, 4.5) {$\Sigma c_{n+2}$} edge [-> , "$2$"] (v)
                                                         edge [->, "$-2$", swap] (w) 
                                                         edge [->, dashed, "$R_{c_{n+2}}$", near end] (c);
            \end{tikzpicture}
        \caption{An example of $\Mm_\varphi$, with the dashed lines representing $u$, and $\C$ is as described in Example \ref{ex:sq2cat}.}
        \label{fig:coneSq2}
    \end{figure}
    It will help to have a diagram for reference, we include a template of this mapping cone in Figure \ref{fig:coneSq2}. 
    
    We begin with objects $c_i\in \Ob(\C)$ with $|c_i|_\C = n$. We define $\Mm_\varphi$ to have the moduli spaces:
    \[\M(\Sigma c_i, e_n) = \begin{cases}
                                P_{c_i} & \mathrm{if }\;\;  \varphi(c_i) =1\\
                                \emptyset &\mathrm{otherwise.}
                            \end{cases}
    \]
    We set $s(P_{c_i}) = 0$. By definition we have then that $u^*(e_n) = \varphi$. In Figure \ref{fig:coneSq2} this is represented by the dashed line from the $\Z/2\Z$-cocycle $\Sigma c_n$ to $e_n$.

    Next, we look at the objects $c_j$ with $|c_j|_\C = n+1$. For each $c_j$, define $M_j^+\subset M$, the set of matched pairs $(P_{j}, P_{j}')$ in $M$ with $\mathrm{source}(P_{j})=c_j$ and $s(P_j) = s(P_j')$. Then we have the following 0-dimensional moduli spaces:
    \[
    \M(\Sigma c_j, e_{n+1}) = \begin{cases}
                            \bigsqcup Q_{P_j}    &\text{where each $Q_{P_j}$ is associated to a} \\
                            &\text{pair $(P_{j}, P_{j}') \in M_j$}\\
                            \emptyset &\text{otherwise.}
                        \end{cases}
    \]
    We set $s(Q_{P_j}) = s(P_{j})$. In Figure \ref{fig:coneSq2}, we have $\sq^1_{c_n}$ acts non-trivially on $c_{n+1}^0$, so we add a dashed line from $\Sigma c_{n+1}^0$ to $e_{n+1}$, and likewise for $c_{n+1}^1$.

    Before we move onto the 0-dimensional moduli spaces $\M(\Sigma c_k,e_{n+2})$ for $|c_k|_\C = n+2$, we will describe the 1-dimensional moduli spaces $\M(\Sigma c_j, e_n)$. We have two types of intervals to describe.

    Let $(P,P') \in M$ with $s(P)\neq s(P')$, with $\mathrm{source}(P)= \mathrm{source}(P') = c_j$ and $\mathrm{target}(P)=c_1, \mathrm{target}(P') =c_2$. Then let there be an interval $I_{(P,P')}$ with $\partial I_{(P,P')} = (P_{c_1}, P)\sqcup (P_{c_2},P')$. evidently the sign assignment for this interval is satisfied as $s(P)\neq s(P')$. We set $f(I_{(P,P')} ) = 0$

    When we have $(P,P')$ with $c_j,c_1,c_2$ as before, but now with $s(P)=s(P')$, we instead add two intervals $I_P, I_{P'}$. These intervals will have boundaries $\partial I_P = (P_{c_1}, P)\sqcup (Q,Q_P)$ and $\partial I_{P'} = (P_{c_2}, P') \sqcup (Q', Q_{P})$. We set $f(I_P) =0$ and $f(I_{P'})=1$. Since $s(Q_P)=s(P)$, these intervals also form anti-commuting squares as $s(\Sigma P) = s(P)+1$, hence $s$ is a sign assignment.

    We have $|\M(\Sigma c_j, e_{n+1})| = |M_j^+|$, hence $u^*(e_{n+1})= \sq^1(\varphi) $.

    Let $c_k$ be an object in degree $n+2$ in $\C$. Since we have non-empty 0-dimensional moduli spaces $\M(\Sigma c_j,e_{n+1})$, we need to assign 1-dimensional moduli spaces to $\M(\Sigma c_k, e_{n+1})$. The first Steenrod square $\Sq^1$, arises as the Bockstein homomorphism of the short exact sequence
    \[
    0\longrightarrow \Z/2\Z \longrightarrow \Z/4\Z \longrightarrow \Z/2\Z \longrightarrow 0;
    \]  
    however, this operation induced on cohomology in coefficients $\Z/2\Z$ can arise as the Bockstein of the sequence
    \[
    0\longrightarrow \Z \longrightarrow \Z \longrightarrow \Z/2\Z\longrightarrow 0,
    \]
    when composed with a change of coefficients, see Definition \ref{def:sq1} for the combinatorial equivalents. In other words we can lift $\sq^1(\varphi)$ to an integral cocycle $\widetilde{\sq}^1$ and then pair up the moduli spaces $\M(c_k,\sq^1(\varphi))$ with a combinatorial boundary matching such that every pair $(R,R')$ has $s(R)\neq s(R')$. Call this set of pairings $M_{\sq^1(\varphi)}$, and then we decide the intervals in $\M(\Sigma c_k, e_{n+1})$ based on the pairs $(R,R') \in M_{\sq^1(\varphi)}$. For each such pair, we have an interval $I_{(R,R')}$ which is defined by $\partial I_{(R,R')} = (Q_{P_1}, R)\sqcup (Q_{P_2}, R')  $. We frame these intervals $f(I_{(R,R')}) =0$.

     Since $f(S)=0$ in $\mathcal{K}_n$, we can always ensure that the compatibility condition determined by $\partial M(\Sigma c_k, e_n)$ is satisfied by simply adding in points $R_{c_k}$ to the moduli spaces $\M(\Sigma c_k, e_{n+2})$.

     Crucially, the statement that we can pull back $u$ to get the Steenrod squares follows by naturality of the mapping cones with respect to the second Steenrod square, see Lemma \ref{lemma:sesConnHom}.  In Figure \ref{fig:coneSq2}, this means must have a point $R_{c_{n+2}}$ in $\M(\Sigma c_{n+2}, e_{n+2})$.
     
     We have $\sq^2(e_n) = e_{n+2}$ by construction, and hence $u^*\circ \Sq^2 = \Sq^2 \circ u^*$ gives $\sq^2(\varphi) = u^*(e_{n+2})$. This means that by simply framing this mapping cone, we are done, and necessarily so, $u^*(e_{n+2})$ will be a cocycle. This fact lets us deal with the objects in degree $n+3$ as well. 

     Since we now have 0-dimensional moduli spaces $\M(\Sigma c_k, e_{n+2})$, we must consider the 0-dimensional moduli spaces $\M(\Sigma c_l, e_{n+3}')$ and the 1-dimensional moduli spaces $\M(\Sigma c_l,e_{n+2})$, for $c_l$ in degree $n+3$. Since $u^*(e_{n+2})$ is a cocycle, we can apply another combinatorial boundary matching and construct these points and intervals exactly as we did for $\M(\Sigma c_j, e_n)$, and, since $e_{n+2}$ has no 0-dimensional moduli spaces with it as the source (and since we are working with the $(n+3)$ truncation of $\C$), we have no compatibility condition to verify, so the framing is inconsequential. We are left now with compatibility conditions corresponding to $\partial \M(\Sigma c_l, e_{n+1})$, but, we can always ensure they are satisfied for the same reason as those corresponding to $\partial \M(\Sigma c_k, e_n)$ were -- we simply add in a point into the moduli space $\M(\Sigma c_l, e_{n+3})$ whenever the condition is not satisfied.

    Hence, the mapping cone can be framed.
\end{proof}

\begin{remark}
    It is worth noting that the $(n+4)$-truncation or any higher degree 1-flow category for $K(\Z/2\Z,n)$ is not just unnecessary for us-- we would also have in general no such mapping cone characterising the Steenrod squares. As stated in \cite{d33}, not every framed 1-flow category can be upgraded to a general framed flow category -- for example, we can have a non-trivial $\Sq^2(\Sq^2(c))$ without requiring that $\Sq^1(\Sq^2(\Sq^1(c)))$ is even non-trivial, let alone equivalent -- and hence no such mapping cone can exist in general.
\end{remark}

% \begin{remark}
%     We also comment on the uniqueness of this construction and its comparison with the special graph structures. Though this is not necessary for our claim, one can see the formula for the second Steenrod squares appearing in the compatibility condition described by $\partial\M(\Sigma c_{n+2}, e_n)$. For instance, by framing the pairs of intervals which are described by a pair $(P,P')$ so that $f(I_P) =0$ and $f(I_{P'})=1$, we effectively account for the $+1$'s corresponding to $D(C)$ in the formula. One could then uniquely define the moduli spaces $\M(\Sigma c_{n+2},e_{n+2})$ by including a 0-framed point whenever a component in the compatibility condition sums to 1, but this does not change things on cohomology. 
% \end{remark}
% 
\section{Spanier-Whitehead Duality}
In this section we define a category $\C^\vee$ in which the objects and moduli spaces are effectively `flipped upside-down', and show that the flipped $0$- and $1$-dimensional moduli spaces can still be framed with some suitable sign assignment and framing.

\begin{definition}
\label{def:dual}
Let $(C,s,f)$ be a framed 1-flow category and let $x,y \in \Ob(\C)$. $\C^\vee$ is defined by:
\begin{itemize}
\item $\Ob(\C^\vee) = \Ob(\C)$, and furthermore $\norm{x^\vee}_{\C^\vee} = -\norm{x}_\C$,
\item $\M(y^\vee,x^\vee)= \M(x,y)$.
\end{itemize}
We also define the sign assignment $s^\vee$. For 0-dimensional $P\in\M(a,b)$:
\[s^\vee (P^\vee) := s(P)\]
And we define the framing $f^\vee$. For circle components $ S \subseteq \morC{a}{c}$,  
        $$f^\vee(S^\vee)=f(S),$$
and for interval components $I \subseteq \morC{a}{c}$ where the boundary is \newline $\partial I = (P,P')\sqcup (Q,Q')$, and where, without loss of generality, we may assume $s(P)=s(P')$, we define 
\[f^\vee(I^\vee) = f(I) + s(P).\]
\end{definition}
We will use the following to prove that this is indeed a framed 1-flow category.
\begin{definition}
First, for $\norm{a} = \norm{d}+3$, let $\partial \morC{a,d}$ be the 1-dimensional framed boundary in $\C$ as in the compatibility condition.
     Let $V$ be the set of \textit{vertices} in $\partial \morC{a}{d}$, more precisely, for $\norm{a}=\norm{b}+1=\norm{c}+2=\norm{d}+3$, we define 
     $$V:=\{ (C,B,A) \; \vert \; C \in \morC{c}{d}, B\in \morC{b}{c}, C\in \morC{a}{b}\} .$$
     Let $V_I , V_J$ be (not necessarily disjoint) subsets of $V$ defined by 
     \begin{align*}
     V_I &:= \{(C,B,A) \;\vert \; s(B)=s(C)\} \;\mathrm{and, }\\
     V_J &:= \{ (C,B,A) \;\vert\; s(B) = s(A)\},
     \end{align*}
    where $s$ is the sign assignment on these 0-dimensional moduli spaces.
\end{definition}
\begin{lemma}
    \label{vertex-triples}
        The sets $V_I, V_J$ are in one-to-one correspondence with intervals of the form $I \times \{A\}$ and $\{C\} \times J$ in $\partial \morC{a}{d}$.
    \end{lemma}
    \begin{proof}
       We remark that for any interval $L$ and $(P,Q) \cup (R,S) = \partial L $, we have, by the requirement on the sign assignment $s(P)+s(Q)+s(R)+s(S)=1$ that exactly one of $(P,Q)$, $(R,S)$ has equal sign assignments. Hence, we have our result almost immediately -- each vertex in $V_I$ touches exactly one edge in $\partial_1\morC{a}{d}$ and each vertex in $V_J$ touches one edge in $\partial_2\morC{a}{d}$ which correspond to intervals of the form $I \times \{C\}$ and $\{A\} \times J$ respectively.
    \end{proof}
\begin{theorem}
    Let $(\C,s,f)$ be a framed 1-flow category. Then $s^\vee$ and $f^\vee$ make $\C^\vee$ into a framed 1-flow category.
\end{theorem}

\begin{proof} 
    It is clear that the sign assignment will still function as required. We focus then on the framing of the intervals.
    
    The compatibility condition on $f^\vee$ will be demonstrated per component. For $S^1$ components, the framing is unchanged, and so their contribution to the overall sum is unchanged. We then have 3 cases to examine which we will define by examining sign assignments at the ``vertices" of the 2n-gon which arises as the boundary of any 2-dimensional moduli space, as in the compatibility condition.
    
    Let $\partial \morC{a,d}$ be a 1-dimensional framed boundary in $\C$, that is, $|a| = |d|+3$. We begin with our sets $V_I$ and $V_J$ as defined in Lemma \ref{vertex-triples} and we use this to define our cases.

In the case of $(C,B,A) \in V_I$, and $s(B) \neq s(A)$, we have to track the contribution of the corresponding interval $I$ to the compatibility condition.\\
    In $\C$, we have that the compatibility framing $\widetilde{f}(I\times\{A\})=f(I)$, however, in $\C'$ this interval is of the form $\{A'\}\times I'$, and hence in the compatibility condition it will contribute:
    \begin{align*}
        \widetilde{f^\vee}(\{A^\vee\}\times I^\vee) &= 
        1+s^\vee(A^\vee) + f^\vee(I^\vee) \\
       &= 1+s(A) + f(I) + s(B) \\
       & = f(I),
    \end{align*}
which follows as it is assumed that $s(B)\neq s(A)$.

In the case of $(C,B,A) \in V_J$ and $s(B)\neq s(C)$, the solution is similar:
    \begin{align*}
    \widetilde{f^\vee}(J^\vee\times \{C^\vee\}) &= f^\vee(J^\vee) \\
       &=f(J) + s(B) \\
       & = f(J) + 1 + s(C),
    \end{align*}
which is of course the contribution of an interval of the form $\{C\}\times J$ to the compatibility condition in $\C$.

Finally, we have to deal with the intersection of the two. In this case, the vertex $(C,B,A)$ with $s(C)=s(B)=s(A)$ corresponds to one interval of the form $I \times \{A\}$ and one interval of the form $\{C\}\times J$. Hence, keeping track of the contribution of \textit{both} intervals, we see that they pair up to produce the desired result:
    \begin{align*}
        \widetilde{f^\vee}(\{A\}\times I) + \widetilde{f^\vee}(J\times \{C\}) &
        = 1+s^\vee(A) + f^\vee(I)+f^\vee(J) \\
       &=1+s(A) + f(I) + s(B)+ f(J) + s(B) \\
       & = f(I) +  1 + s(C) + f(J),
    \end{align*}
where we make use of the equivalence in the sign assignments in the last step. 

So, the framing $f^\vee$ contributes the same values to the compatibility condition as $f$ when paired with the sign assignment $s^\vee(P^\vee) = s(P)$, and hence $(\C^\vee, s^\vee, f^\vee)$ is framed.
\end{proof}
The following theorem parallels Switzer's treatment of the SW, or S-dual, in \cite[Ch. 14]{Switzer1975-pm}. The existence of a map is instead resolved as the existence of a framed mapping cone 1-flow category.

\begin{theorem} \label{thm: swDual}
    There exists a `mapping cone' framed 1-flow category $\Mm_\mu$ which fits into the short exact sequence of 1-flow categories:
        \[
        0\longrightarrow \oldcal{S}^0 \longrightarrow \Mm_\mu \longrightarrow \Sigma (\C^\vee\otimes \C)\longrightarrow 0
        \]
    $\M_\mu$ determines an isomorphism $H_i(\C^\vee;\Z) \rightarrow H^{-i}(\C;\Z) $. 
    
\end{theorem}
\begin{proof}
We first define $\M_\mu$ and show that it is framed. The objects of $\Mm_\mu$ are comprised of a disjoint union of the objects of $\Sigma(\C^\vee\otimes \C)$, and $\Ob{\oldcal{S^0}} = \{s_0\}$, which is the 1-flow category with a unique object, in degree 0. 

All moduli spaces and their framings in the suspension of $\C^\vee\otimes\C$ are the same in $\Mm_\mu$, but we have additional 0-dimensional moduli spaces 
\[\M(\Sigma(c^\vee,c),s_0) = \{P_c\},
\]
 with $s(P_c) = \frac{1}{2}|c|(|c|+1)$, which we denote $t_{|c|+1}$ (i.e the $(|c|+1)^{\mathrm{th}}$ triangular number). For every point $Q$ in any 0-dimensional moduli space $\M_\C(a,b)$, we have in $\Mm_\mu$ an interval $I_Q \cong \M((b^\vee,a), s_0)$, which we show in Figure \ref{fig:muInterval}, and we set the framing $f(I_Q) = s(Q)$.
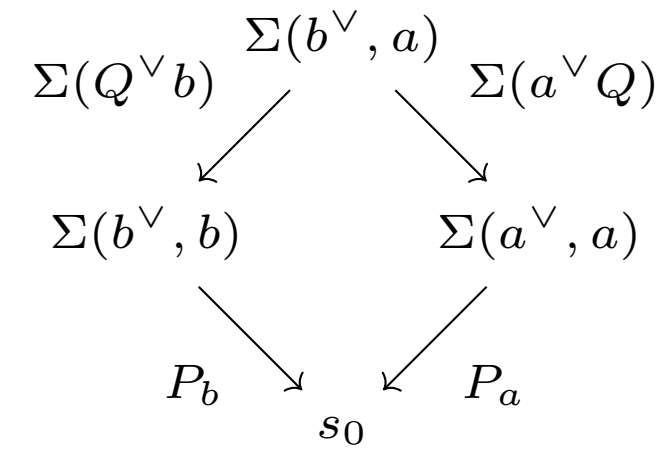
\begin{figure}[htbp]
\centering
    \begin{tikzpicture}[font=\tiny, transform shape, scale=1.3]
    \begin{scope}
        \node (a) at (0,0) {$\Sigma(b^\vee,a)$};%
        \node (b) at (-1,-1) {$\Sigma(b^\vee,b)$} edge ["$\Sigma (Q^\vee b)$" , <-](a) ;%
        \node (c) at (1,-1) {$\Sigma(a^\vee,a)$} edge [ swap,"$\Sigma (a^\vee Q)$", <-] (a);%
        \node (d) at (0,-2) {$s_0$} edge ["$P_{b}$", <-] (b)
                                        edge [ swap,"$P_{a}$", <-] (c);
    \end{scope}

    \end{tikzpicture}
    \caption{The subcategories corresponding to intervals $I_Q$ in $\Mm_\mu$.}
    \label{fig:muInterval}
\end{figure}
For the sign assignments, we need 
\[
s(P_{b}) + s(P_{a}) + s(\Sigma(Q^\vee b)) + s(\Sigma(a^\vee Q)) = 1.
\] %
The latter two summands are defined using the product and suspensions of framed 1-flow categories and give %
\[s(\Sigma(Q^\vee b) + s(\Sigma(a^\vee Q) = (1+s(Q)) + (1 + (|a| + s(Q)))  = |a| .
\]%
Since the defined sign assignments on $P_{b}$ and $P_{a}$ sum to $|a|+1$, $s$ is a sign assignment.

For every interval $J\subseteq \M(a,c) $ in $\C$ with boundary $(Q_1,Q_3)\sqcup (Q_2,Q_4)$ (where the pair with equal sign assignment are denoted by $(Q_1,Q_3)$) there will be a compatibility condition to satisfy in $\partial \M(\Sigma(a^\vee,c),s_0)$. This yields a subcategory of $M_\mu$, shown in Figure \ref{fig:muComp}, and in the boundary of this moduli space is an octagon, as in Figure \ref{fig:muComp2}. To make these diagrams easier to read the `$\Sigma$'s have suppressed.
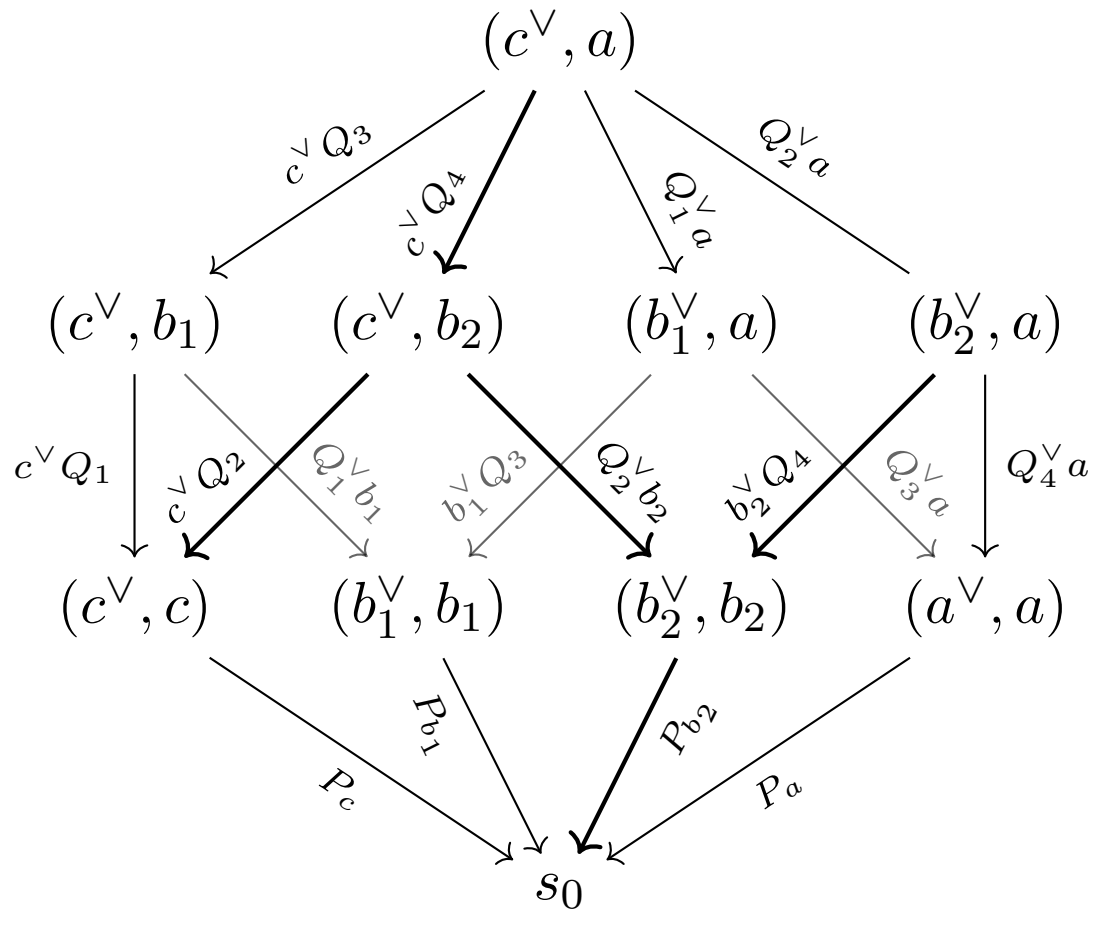
\begin{figure}[htbp]
    \centering
    \begin{tikzpicture}[every edge quotes/.style = {auto, font=\ttfamily\tiny,
  text=black, fill=none}, scale =0.75]
        \begin{scope}[scale = 2.5]
            \node (s0) at (0,0) {$s_0$};
            \node (c0c0) at (-1.5,1) {$(c^\vee,c)$} edge ["$P_{c}$",swap,->,sloped] (s0);
            \node (c11c11) at (-0.5,1) {$(b_1^\vee,b_1)$} edge ["$P_{b_1}$", near start,->,sloped,swap] (s0);
            \node (c12c12) at (0.5,1) {$(b_2^\vee,b_2)$} edge [thick, near start,"$P_{b_2}$",->,sloped, near start,swap] (s0) ;
            \node (c2c2) at (1.5,1) {$(a^\vee,a)$} edge ["$P_{a}$",->,sloped,swap] (s0);
            \node (c11c0) at (-1.5,2) {$(c^\vee,b_1)$} edge [swap, "$c^\vee Q_1$", ->] (c0c0)
                                                         edge [near end,sloped, "$Q_1^\vee b_1$", opacity = 0.6, ->] (c11c11);
            \node (c12c0) at (-0.5,2) {$(c^\vee,b_2)$} edge [near end, "$c^\vee Q_2$", thick,sloped, ->] (c0c0)
                                                         edge [near end, "$Q_2^\vee b_2$", thick, ->,sloped] (c12c12);
            \node (c2c11) at (0.5,2) {$(b_1^\vee,a)$}  edge [sloped, near end, "$b_1 ^\vee Q_3$",opacity = 0.6,->] (c11c11)
                                                         edge [near end, sloped, "$Q_3^\vee a$", opacity = 0.6, ->] (c2c2);
            \node (c2c12) at (1.5, 2) {$(b_2^\vee,a)$}  edge [near end,sloped, "$b_2^\vee Q_4$", thick, ->] (c12c12)
                                                       edge [ "$Q_4^\vee a$", ->] (c2c2);
            
            \node (c2c0) at (0,3) {$(c^\vee,a)$} edge [sloped, "$c^\vee Q_3$", ->] (c11c0)
                                                  edge [thick, near end,sloped, "$c^\vee Q_4$", ->] (c12c0)
                                                  edge [near end, "$Q_1^\vee a$", ->,sloped] (c2c11)
                                                  edge ["$Q_2^\vee a$"->,sloped] (c2c12);
        \end{scope}
        
    \end{tikzpicture}
    \caption{The 1-flow subcategory $\M(\Sigma(c^\vee,a),s_0)$}
    \label{fig:muComp}
\end{figure}
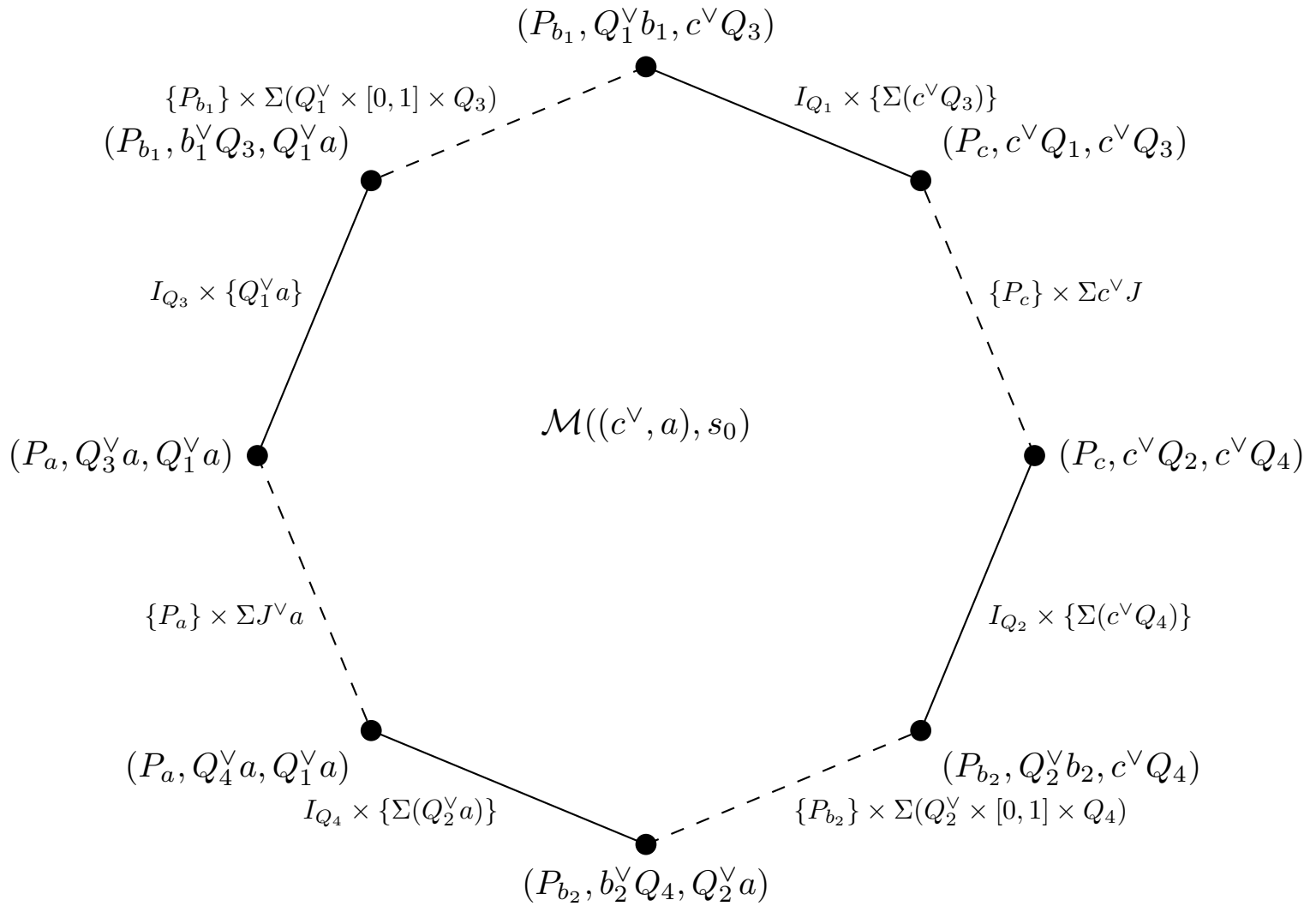
\begin{figure}[htbp]
    \centering
    \begin{tikzpicture}[scale = 0.75, transform shape]
        \begin{scope}[font =\small]
            \newdimen\R
            \R=4cm
            \node[anchor = south] {$\M ((c^\vee,a),s_0)$};
            \foreach \x/\l/\d in
                { 45/{$\{P_c\}\times \Sigma c^\vee J$}/dashed,
                90/{$I_{Q_1}\times \{\Sigma (c^\vee Q_3)\}$}/,
                135/{$\{P_{b_1}\}\times \Sigma(Q_1^\vee\times [0,1]\times Q_3)$}/dashed,
                180/{$I_{Q_3}\times \{Q_1^\vee a\}$}/,
                225/{$\{P_a\}\times \Sigma J^\vee a$}/dashed,
                270/{$I_{Q_4}\times \{\Sigma(Q_2^\vee a)\}$}/,
                315/{$\{P_{b_2}\}\times \Sigma( Q_2^\vee\times[0,1]\times Q_4)$}/dashed,
                360/{$I_{Q_2}\times \{\Sigma(c^\vee Q_4)\}$}/%
                }
            \draw[postaction={decorate}] ({\x-45}:\R) edge[-,\d] node[auto,swap, scale = 0.8]{\l} (\x:\R);
            \foreach \x/\l/\p in
                { 45/{$(P_c,c^\vee Q_1,c^\vee Q_3)$}/above right,
                90/{$(P_{b_1},Q_1^\vee b_1,c^\vee Q_3)$}/above,
                135/{$(P_{b_1}, b_1^\vee Q_3, Q_1^\vee a)$}/above left,
                180/{$(P_a, Q_3^\vee a,Q_1^\vee a)$}/left,
                225/{$(P_a, Q_4^\vee a, Q_1^\vee a)$}/below left,
                270/{$(P_{b_2},b_2^\vee Q_4,Q_2^\vee a)$}/below,
                315/{$(P_{b_2},Q_2^\vee b_2,c^\vee Q_4)$}/below right,
                360/{$(P_c,c^\vee Q_2,c^\vee Q_4)$}/right
                }
            \node[inner sep=2pt,circle,draw,fill,label={\p:\l}] at (\x:\R) {};
        \end{scope}
    \end{tikzpicture}
    \caption{The octagon $\partial \M((c^\vee,a),s_0)$ with intervals labelled}
    \label{fig:muComp2}
\end{figure}
 The compatibility condition corresponding to the octagon $C=\partial \M((c^\vee,a),s_0)$ yields the sum:
\[
\Tilde{f}(K) = \sum_{I\times\{R\}}f(I) +  \sum_{\{P\}\times J} 1 + s(P) + f(J).
\]
Of the first summand, we have four intervals $\left\{I_{Q_i}\right\}_{i=1}^4$, whose framings sum to $1$. Of the second summand, we have $\sum 1 + s(P_c) =1$, and we then have the sum of the framings $\sum f(J)$:
\begin{align*}
    f(\Sigma(J^\vee a)) + f(\Sigma (Q_1^\vee \otimes Q_3)) +f(\Sigma(c^\vee J)) + f(\Sigma(Q_2^\vee \otimes Q_4)).  
\end{align*}
Following our definitions carefully indeed shows that the compatibility condition is satisfied for these components.

Lastly, for circle components $S\subseteq \M(a,c)$, the contributions to the compatibility condition from $S_{c^\vee}$ in $\M((a,c^\vee) , (c,c^\vee))$ will pair up with the dual $aS$ in $\M((a,c^\vee) , (a,a^\vee))$. So, $\Mm_\mu$ is framed.

In order to see that a perfect pairing is produced, we note that the 0-dimensional moduli spaces of $\Mm_\mu$ define a chain map 
\[\mu:C_*(\C\otimes \C^\vee;\Z) \to \Z,
\]
 where $\mu$ is defined from the moduli spaces $\M(\Sigma(\C\otimes\C^\vee),s_0)$. For pairs $\Sigma(c,c^\vee)\in \Ob(\M_\mu)$ this map is given by $\mu((c,c^\vee)) = -1^{t_{|c|+1}}$, which is clearly non-degenerate and identifies the cochain groups of $\C$ with the chain groups of $\C^\vee$, and hence it does so on (co)homology as well.
\end{proof}

Our dual has the following properties:

 \begin{lemma}
     Let $\C,\D$ be framed 1-flow categories, then we have identities $(\C^\vee)^\vee = \C$, and $(\C\otimes\D)^\vee = \C^\vee \otimes \D^\vee$. Furthermore, if we have a short exact sequence of framed 1-flow categories
     \[
     \D \rightarrow \Mm_f \rightarrow \Sigma \C,
     \]
     then we get from the dual a short exact sequence of framed 1-flow categories
     \[
     \Sigma^{-1}\C^\vee \rightarrow \Mm_f^\vee \rightarrow \D^\vee.
     \]
 \end{lemma}

 \begin{proof}
 The proof of the first two statements follows by straightforward checks of the definitions. The last statement follows by dualising the mapping cone $\Mm_f$. There is an obvious identification of the objects in $\Mm_f^\vee$ with $\Ob(\Sigma^{-1} \C^\vee) \sqcup \Ob(\D^\vee)$, and since $\Sigma \C$ is upwards-closed in $\Mm_f$, its dual is downwards-closed in $(\Mm_f)^\vee$, and likewise for the subcategory $\D$ and $\D^\vee$. For an interval $I$ in $\C$ with boundary $(P,P')\sqcup (Q,Q')$ with $s(P)=s(P')$, we have
\begin{align*}
f((\Sigma I)^\vee) &= (s(P) + 1) + (s(P)+ s(Q) + f(I))\\
 &= (s(P) + s(Q) + 1) + (s(P) + f(I)) \\
 &= s(P') (s(P')+s(Q') +f(I))\\
 &= f(\Sigma^{-1}(I^\vee)),
\end{align*}
hence the framings of $\Sigma^{-1}\C^\vee$ are equivalent to their image in $(\Mm_f)^\vee$ as well.
\end{proof}
These identities align with what one expects from the SW-dual and the product, see \cite{Switzer1975-pm} for example. We claim that $\C^\vee$ is the correct equivalent of a Spanier-Whitehead dual to $\C$. This statement itself might take a little more work, in particular, one might need to consider a notion of homotopy equivalence for framed 1-flow categories in order to show something like uniqueness, however for practical purposes, one needs only to know how the dual cooperates with the Steenrod square.
\begin{lemma}
    Let $\mathcal{K}_n$ be the framed 1-flow category described in Definition \ref{def:eilenberg}. Then the following diagram commutes for all $k\in \Z$
    \begin{center}
        \begin{tikzpicture}[ font =\small, transform shape]
            \node (a) at (-2.5,0) {$H^k(\mathcal{K}_n;\mathbb{F}_2)$};
            \node (b) at (2.5,0) {$\Hom (H^{-k}(\mathcal{K}_n^\vee;\mathbb{F}_2),\mathbb{F}_2)$} edge[<-, "$\simeq$"] (a);
            \node (c) at (-2.5,-2) {$H^{k+2}(\mathcal{K}_n;\mathbb{F}_2)$} edge[<-, "$\Sq^2$"] (a) ;
            \node (d) at (2.5,-2) {$\Hom(H^{-k-2}(\mathcal{K}_n^\vee;\mathbb{F}_2),\mathbb{F}_2)$.} edge[<-, "$\simeq$"](c) edge[<-, "$(\Sq^2)^*$",swap] (b);
        \end{tikzpicture}
    \end{center}
\end{lemma}
\begin{proof}
 The framed 1-flow category for $\mathcal{K}_n^\vee$ is represented by the following diagram:
    
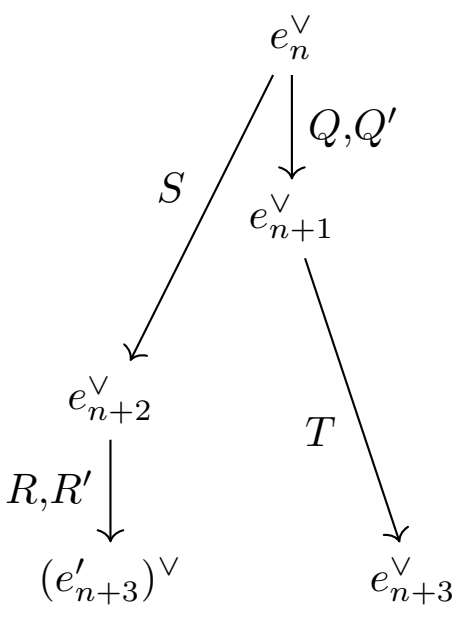
\begin{center}
    \begin{tikzpicture}[font = \small, transform shape, scale =0.8]
        \node (a) at (0,0) {$e_{n}^\vee$};
        \node (b) at (0,-1.5) {$e_{n+1}^\vee$} edge [<-, swap, "$Q{,}Q'$"] (a);
        \node (c) at (-1.5,-3) {$e_{n+2}^\vee$} edge [<-, "$S$"] (a);
        \node (d) at (1,-4.5) {$e_{n+3}^\vee$} edge [<-, "$T$"] (b);
        \node (e) at (-1.5,-4.5) {$(e_{n+3}')^\vee$} edge [<-, "$R{,}R'$"] (c);       
    \end{tikzpicture}
\end{center}
The horizontal arrows can be constructed as a composition of the isomorphism of Theorem \ref{thm: swDual}, and the Universal Coefficients theorem which ensures that, for instance, $H_{-k}(\C^\vee;\mathbb{F}_2)\simeq \Hom(H^{-k}(\C;\mathbb{F}_2),\mathbb{F}_2)$. Looking specifically at $[e_n]\in H^n(\mathcal{K};\mathbb{F}_2)$, we can follow the upper horizontal arrow where $[e_n]$ maps to $\Hom([e_n^\vee]^*,\mathbb{F}_2)$, which in turn maps to $\Hom([e_{n+2}^\vee], \mathbb{F}_2)$, as $\Sq^2([e_{n+2}^\vee]) = [e_n^\vee]$. In the other direction, $\Sq^2([e_n]) = [e_{n+2}]$, which corresponds to $\Hom([e_{n+2}^\vee]^*, \mathbb{F}_2)$. The other generators either have a trivial $\Sq^2$ or commute for similar reasons, and hence the diagram commutes.
 \end{proof}
\begin{proposition}
    Let $\C$ be a framed 1-flow category, then the following diagram commutes:
    \begin{center}
        \begin{tikzpicture}[ font =\small, transform shape]
            \node (a) at (-2.5,0) {$H^k(\C;\mathbb{F}_2)$};
            \node (b) at (2.5,0) {$\Hom (H^{-k}(\C^\vee;\mathbb{F}_2),\mathbb{F}_2)$} edge[<-, "$\simeq$"] (a);
            \node (c) at (-2.5,-2) {$H^{k+2}(\C;\mathbb{F}_2)$} edge[<-, "$\Sq^2$"] (a) ;
            \node (d) at (2.5,-2) {$\Hom(H^{-k-2}(\C^\vee;\mathbb{F}_2),\mathbb{F}_2)$.} edge[<-, "$\simeq$"](c) edge[<-, "$\Hom(\Sq^2)$",swap] (b);
        \end{tikzpicture}
    \end{center}
\end{proposition}
\begin{proof}
 Let  $z\in H^n(\C;\Z/2\Z)$ be a cohomology class represented by a cocycle $\varphi \in C^n(\C;\mathbb{F}_2)$ then we analyse the following diagram which we construct given the cocycle $\varphi$ so that the distal arrows arise from the mapping cone $\M_\varphi$ which arises from Proposition \ref{prop:altDefSq2}, and its dual. 

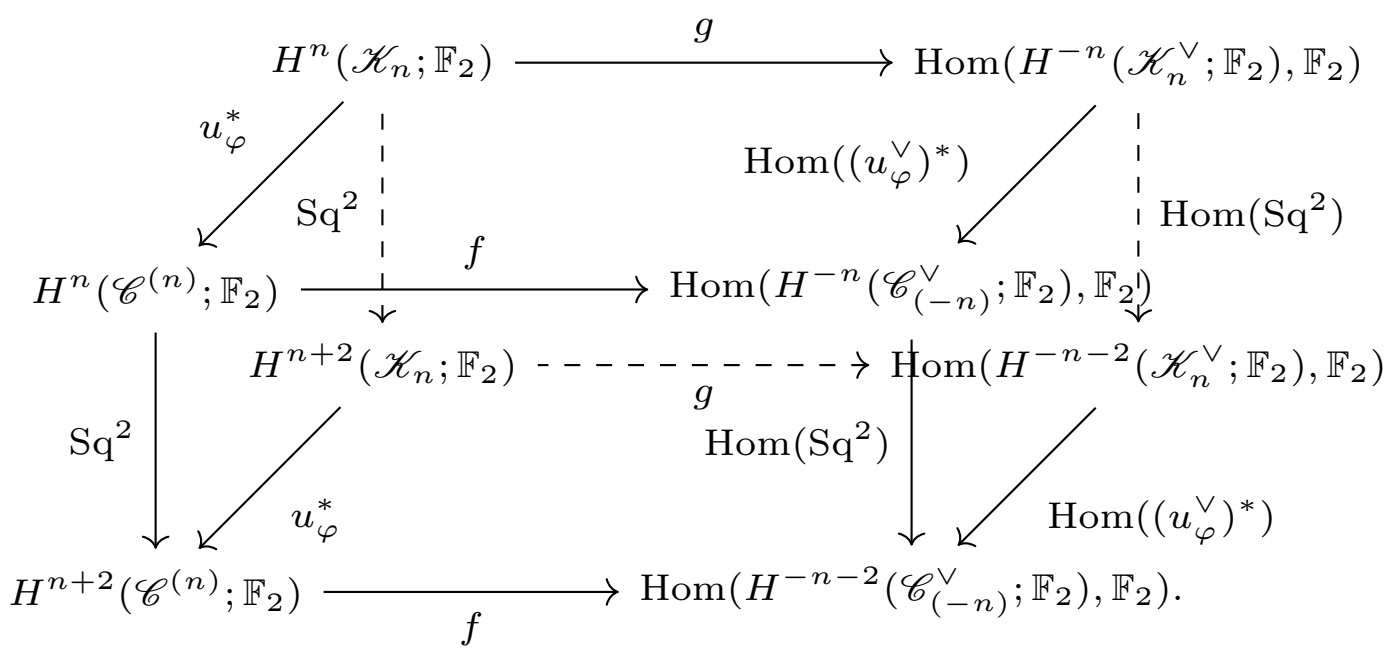
\begin{center}
    \begin{tikzpicture}[font=\tiny, transform shape]
            \node (a1) at (-2.5,0) {$H^n(\C^{(n)};\mathbb{F}_2)$};
            \node (b1) at (2.5,0) {$\Hom (H^{-n}(\C^\vee_{(-n)};\mathbb{F}_2),\mathbb{F}_2)$} edge[<-, "$f$", swap] (a1);
            \node (c1) at (-2.5,-2) {$H^{n+2}(\C^{(n)};\mathbb{F}_2)$} edge[<-, "$\Sq^2$"] (a1) ;
            \node (d1) at (2.5,-2) {$\Hom(H^{-n-2}(\C_{(-n)}^\vee;\mathbb{F}_2),\mathbb{F}_2)$.} edge[<-, "$f$"](c1) edge[<-, "$\Hom(\Sq^2)$"] (b);
            \node[behind path] (a2) at (-1,1.5) {$H^n(\mathcal{K}_n;\mathbb{F}_2)$} edge[->, "$u_\varphi^*$",swap] (a1);
            \node[behind path] (b2) at (4,1.5) {$\Hom (H^{-n}(\mathcal{K}_n^\vee;\mathbb{F}_2),\mathbb{F}_2)$} edge[<-, "$g$", swap, behind path] (a2) edge[->, "$\Hom((u_\varphi^\vee)^*)$",swap, near end](b1);
            \node[behind path] (c2) at (-1,-0.5) {$H^{n+2}(\mathcal{K}_n;\mathbb{F}_2)$} edge[<-, "$\Sq^2$", behind path, dashed] (a2)  edge [->, "$u_\varphi^*$" ] (c1);
            \node[behind path] (d2) at (4,-0.5) {$\Hom(H^{-n-2}(\mathcal{K}_n^\vee;\mathbb{F}_2),\mathbb{F}_2)$} edge[<-, "$g$", behind path ,dashed](c2) edge[<-, "$\Hom(\Sq^2)$",swap, behind path, dashed] (b2) edge [->, "$\Hom((u_\varphi^\vee)^*)$"] (d1);            
    \end{tikzpicture}
\end{center}

We want to show that the front face of the cube commutes. Using $f,g$ to represent the isomorphisms induced by the compositions of the isomorphisms given by the mapping cone $\M_\mu$ and the isomorphisms given by application of the Universal Coefficients theorem, we wish to show that for some cohomology class $[\varphi] \in H^n(\C;\mathbb{F}_2)$, we have that $(\Sq^2)^* f([\varphi]) = f\Sq^2(\varphi)$. We prove this by effectively diagram chasing across the whole cube.

We start by pulling back $[\varphi]$ along $u^*_\varphi$, and then prove this statement from there. The top face commutes by naturality of the Spanier-Whitehead duality and naturality of the isomorphism due to the Universal Coefficients theorem, and hence we can write 
$$\Hom(\Sq^2) fu^*([e_n]) = \Hom(\Sq^2)\Hom( (u_{\varphi}^\vee)^*) g([e_n]).$$
The right face commutes, by naturality of the second Steenrod square with respect to mapping cones, and the back face commutes by the previous Lemma. This leaves us with $\Hom( (u_{\varphi}^\vee)^*) g\Sq^2([e_n])$. Since the bottom face commutes for the same reasons that the top face does, we get to $f u_\varphi^*\Sq^2([e_n])$, and then Proposition \ref{prop:altDefSq2} means that we know $u_\varphi^*\Sq^2([e_n]) = \Sq^2([\varphi])$, hence we have the desired commutativity on the front face for our choice of cohomology class.

It follows then by inclusion of $\C^{(n)}$ into $\C$, that the diagram commutes for this choice of $z$, and since the choice of cohomology class was arbitrary, the front face commutes generally also. 
 \end{proof}
\section{Discussion}
 Missing from this work are discussions of uniqueness. We believe that a correct notion of homotopy equivalence can be constructed in the context of framed 1-flow categories, and in some homotopy category of framed 1-flow categories, one should find that the product is unique, and moreover commutative as well, and that the dual is unique. One possible notion of uniqueness is by  considering these mapping cones as identifying two framed 1-flow categories if and only if the mapping cone reduces to a contractible chain complex. This should force, by naturality, that the framed 1-flow categories on either side of this mapping cone have effectively the same cohomology and Steenrod squares. These details appear in \cite{thesis}, where we will explore this idea further.

 Framed flow categories have found use in other contexts,\cite{manolescu2021knotfloerstablehomotopy, bhattacharyya2025lipshitzsarkarstablehomotopytype}, the former of which has more recently seen work on the construction of a Steenrod square for link Floer homology via framed 1-flow categories in \cite{tao2026steenrodsquarelinkfloer}. This was also used by Sch\"utz in \cite{schuetz2022secondsteenrodsquaresodd} to construct two second Steenrod squares on odd Khovanov homology, which has since been shown to coincide with the second Steenrod square on the odd spectrum by Rajapakse \cite{rajapakse2026steenrodsquaresoddkhovanov}. In upcoming work \cite{schutzJohnson} we will use the 1-flow category dual to resolve a conjecture of Sch\"utz about the nature of the two second Steenrod squares on odd Khovanov homology when applied to the mirror of a knot, and we hope that expanding the tools available in this context allow for similar advances where framed 1-flow categories arise.
\bibliographystyle{amsalpha}
\bibliography{references}

\end{document}